\documentclass{article}

\usepackage{arxiv}

\usepackage[utf8]{inputenc} 
\usepackage[T1]{fontenc}    
\usepackage{hyperref}       
\usepackage{url}            
\usepackage{booktabs}       
\usepackage{amsfonts}       
\usepackage{nicefrac}       
\usepackage{microtype}      
\usepackage{lipsum}		
\usepackage{graphicx}
\usepackage{natbib}
\usepackage{doi}
\usepackage{enumerate}
\usepackage{enumitem}
\usepackage{amsmath, amsthm, amssymb, algorithm, algpseudocode}
\usepackage[table]{xcolor}
\usepackage{tikz}
\usepackage{bm}

\renewcommand{\to}{\longrightarrow}

\newtheorem{theorem}{Theorem}[section]

\newtheorem{proposition}[theorem]{Proposition}

\newtheorem{definition}[theorem]{Definition}
\newtheorem{example}[theorem]{Example}

\newtheorem{problem}[theorem]{Problem}
\newtheorem{application}{Application}
\newtheorem{remark}[theorem]{Remark}

\newcommand\mystyle{\everymath{\displaystyle}}
\mystyle
\title{Complete Structural Analysis of $q$-Heisenberg Algebras: \\
       Homology, Rigidity, Automorphisms, and Deformations}

\author{\href{https://orcid.org/0000-0002-3816-5287}{\includegraphics[scale=0.06]{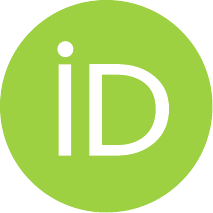}\hspace{1mm}M.H.M.~Rashid}\thanks{Corresponding Author} \\
	Department of Mathematics\&Statistics\\Faculty of Science P.O.Box(7)\\
	Mutah University University\\
	Mutah-Jordan \\
	\texttt{mrash@mutah.edu.jo}
}

\renewcommand{\shorttitle}{\textit{arXiv} Template}

\hypersetup{
pdftitle={Complete Structural Analysis of $q$-Heisenberg Algebras: \\
       Homology, Rigidity, Automorphisms, and Deformations},
pdfsubject={q-bio.NC, q-bio.QM},
pdfauthor={M.H.M.Rashid},
pdfkeywords={q-Heisenberg algebra, iterated Ore extension, global dimension,
          PBW basis, deformation quantization, graded automorphisms,
          Hilbert series, noncommutative algebra}}

\begin{document}
\maketitle

\begin{abstract}
	This paper establishes several fundamental structural properties of the $q$-Heisenberg algebra $\mathfrak{h}_n(q)$, a quantum deformation of the classical Heisenberg algebra. We first prove that when $q$ is not a root of unity, the global homological dimension of $\mathfrak{h}_n(q)$ is exactly $3n$, while it becomes infinite when $q$ is a root of unity. We then demonstrate the rigidity of its iterated Ore extension structure, showing that any such presentation is essentially unique up to permutation and scaling of variables. The graded automorphism group is completely determined to be isomorphic to $(\mathbb{C}^*)^{2n} \rtimes S_n$. Furthermore, $\mathfrak{h}_n(q)$ is shown to possess a universal deformation property as the canonical PBW-preserving deformation of the classical Heisenberg algebra $\mathfrak{h}_n(1)$. We compute its Hilbert series as $(1-t)^{-3n}$, confirming polynomial growth of degree $3n$, and establish that its Gelfand--Kirillov dimension coincides with its classical Krull dimension. These results are extended to a generalized multi-parameter version $\mathfrak{H}_n(\mathbf{Q})$, and illustrated through detailed examples and applications in representation theory and deformation quantization.
\end{abstract}

\section{Introduction}

The Heisenberg algebra, originating from Werner Heisenberg's pioneering work in quantum mechanics \cite{Heisenberg1925}, has become a cornerstone of both Lie theory and noncommutative algebra. Its universal enveloping algebra---the classical Heisenberg algebra $\mathfrak{h}_n(1)$---exhibits a rich structure that has been extensively studied from various perspectives, including automorphism groups \cite{Zhang2003, Zhou2021}, Rota-Baxter operators \cite{Ji2019}, and representation theory. In recent decades, quantum deformations of classical algebras have emerged as central objects in mathematical physics and noncommutative geometry, leading to the natural $q$-deformation known as the $q$-Heisenberg algebra $\mathfrak{h}_n(q)$ \cite{Wess2000, Rosenberg1995}.

\subsection{Scholars' Contributions}
The study of $q$-Heisenberg algebras builds upon foundational work in several areas. Berger \cite{Berger1992} established early quantum versions of the Poincaré--Birkhoff--Witt theorem, while Rosenberg \cite{Rosenberg1995} explored their geometric significance in noncommutative algebraic geometry. Zhang \cite{Zhang2024} recently provided a comprehensive analysis of Gröbner--Shirshov bases for $\mathfrak{h}_n(q)$, establishing its basic structural properties and PBW basis. The general theory of solvable polynomial algebras and their Gr\"{o}bner bases, developed by Li \cite{Li2011, Li2021, Li2002, Li2018}, provides essential computational tools. Homological aspects draw upon the work of McConnell and Robson \cite{McConnell2001}, Goodearl and Warfield \cite{Goodearl2004}, and the regularity theories of Artin--Schelter \cite{Artin1990}, Stafford \cite{Stafford1994}, and Yekutieli--Zhang \cite{Yekutieli2020}.

\subsection{Core of Study}
This paper investigates deep structural properties of the $q$-Heisenberg algebra $\mathfrak{h}_n(q)$ that go beyond previously established results. We focus on five fundamental aspects:
\begin{enumerate}
    \item The precise determination of global homological dimension (Theorem \ref{THEOREM-A1}), distinguishing between the generic case ($q$ not a root of unity) and the root-of-unity case.
    \item The rigidity of its iterated skew-polynomial structure (Theorem \ref{THEOREM-B}), showing that any presentation as an iterated Ore extension is essentially equivalent to the canonical one.
    \item The complete description of its graded automorphism group (Theorem \ref{THEOREM-C}).
    \item Its universal property as a deformation of the classical Heisenberg algebra (Theorem \ref{THEOREM-D}).
    \item Its Hilbert series and growth behavior (Theorem \ref{THEOREM-E}), confirming polynomial growth of degree $3n$.
\end{enumerate}

\subsection{Applications}
The structural results obtained have significant applications in several domains:
\begin{itemize}
    \item \textbf{Representation Theory:} The rigidity theorem enables systematic classification of module families, particularly weight modules with one-dimensional weight spaces.
    \item \textbf{Deformation Quantization:} The universal deformation property establishes $\mathfrak{h}_n(q)$ as the natural quantization of classical coupled harmonic oscillator systems.
    \item \textbf{Computational Algebra:} The PBW basis and Ore extension structure facilitate Gröbner basis computations and algorithmic solutions to ideal membership problems.
    \item \textbf{Noncommutative Geometry:} The homological regularity properties (Auslander-regular, Cohen--Macaulay) position $\mathfrak{h}_n(q)$ as a model for noncommutative affine spaces.
\end{itemize}
Section \ref{Sec:GISPS} presents a generalized multi-parameter version that extends these results to a broader class of algebras, while Section 6 provides explicit numerical implementations and applications.

\subsection{Significance of Study}
The $q$-Heisenberg algebra serves as a paradigmatic example bridging several mathematical disciplines. Its study reveals:
\begin{itemize}
    \item How deformation parameters affect homological invariants, with a sharp dichotomy between the generic and root-of-unity cases.
    \item The remarkable rigidity of Ore extension structures in certain noncommutative domains.
    \item The interplay between combinatorial properties (PBW bases, Hilbert series) and algebraic structures (global dimension, automorphisms).
    \item The universality of $q$-deformations in preserving PBW structures during quantization.
\end{itemize}
These insights contribute to the broader understanding of noncommutative algebras, quantum groups, and deformation quantization.

\subsection{Organization of the Paper}
The paper is organized as follows. Section 2 recalls preliminary material on solvable polynomial algebras, $q$-Heisenberg algebras, iterated Ore extensions, and homological regularity properties. Section 3 presents a generalized multi-parameter version and establishes its iterated skew-polynomial structure. Section 4 contains the main structural theorems for $\mathfrak{h}_n(q)$. Section 5 provides a detailed example for $n=2$ illustrating all theoretical results. Section 6 develops applications to module classification and deformation quantization, including concrete numerical computations. Section 7 concludes with open problems and directions for further research.

The results presented herein significantly advance our understanding of $q$-Heisenberg algebras, providing both theoretical depth and practical tools for researchers in noncommutative algebra, representation theory, and mathematical physics.
\section{Preliminaries}
\label{sec:preliminaries}

This section recalls the essential definitions, notation, and foundational results that will be used throughout the paper. We focus on the concepts of solvable polynomial algebras, \(q\)-Heisenberg algebras, and their fundamental properties as established in the literature \cite{Li2011, Zhang2024, Goodearl2004, McConnell2001}.

\subsection{Solvable Polynomial Algebras and PBW Bases}

Let \(K\) be a field and let \(A = K\langle a_1, \dots, a_n \rangle\) be a finitely generated associative \(K\)-algebra. Suppose \(A\) has a \emph{Poincar\'e--Birkhoff--Witt (PBW)} \(K\)-basis
\[
\mathcal{B} = \{ a^{\beta} = a_1^{\beta_1} \cdots a_n^{\beta_n} \mid \beta = (\beta_1, \dots, \beta_n) \in \mathbb{N}^n \}.
\]
Let \(\prec\) be a total order on \(\mathcal{B}\). For a nonzero element \(f \in A\), written uniquely as
\[
f = \lambda_1 a^{\beta^{(1)}} + \lambda_2 a^{\beta^{(2)}} + \cdots + \lambda_m a^{\beta^{(m)}},
\]
with \(a^{\beta^{(1)}} \prec a^{\beta^{(2)}} \prec \cdots \prec a^{\beta^{(m)}}\), \(\lambda_j \in K^*\), we denote by \(\mathbf{LM}(f) = a^{\beta^{(m)}}\) the \emph{leading monomial}, by \(\mathbf{LC}(f) = \lambda_m\) the \emph{leading coefficient}, and by \(\mathbf{LT}(f) = \lambda_m a^{\beta^{(m)}}\) the \emph{leading term}.

\begin{definition}[Monomial order \cite{Li2011}]
A total order \(\prec\) on the PBW basis \(\mathcal{B}\) of \(A\) is called a \emph{monomial order} if it satisfies:
\begin{enumerate}
    \item \(\prec\) is a well-order (every nonempty subset of \(\mathcal{B}\) has a minimal element);
    \item For all \(a^{\gamma}, a^{\alpha}, a^{\beta}, a^{\eta} \in \mathcal{B}\), if \(a^{\gamma} \neq 1\), \(a^{\beta} \neq a^{\gamma}\), and \(a^{\gamma} = \mathbf{LM}(a^{\alpha} a^{\beta} a^{\eta})\), then \(a^{\beta} \prec a^{\gamma}\) (hence \(1 \prec a^{\gamma}\) for all \(a^{\gamma} \neq 1\));
    \item For \(a^{\gamma}, a^{\alpha}, a^{\beta}, a^{\eta} \in \mathcal{B}\), if \(a^{\alpha} \prec a^{\beta}\), \(\mathbf{LM}(a^{\gamma} a^{\alpha} a^{\eta}) \neq 0\), and \(\mathbf{LM}(a^{\gamma} a^{\beta} a^{\eta}) \notin \{0,1\}\), then
    \[
    \mathbf{LM}(a^{\gamma} a^{\alpha} a^{\eta}) \prec \mathbf{LM}(a^{\gamma} a^{\beta} a^{\eta}).
    \]
\end{enumerate}
\end{definition}

\begin{definition}[Solvable polynomial algebra \cite{Li2011}]
A finitely generated \(K\)-algebra \(A = K\langle a_1, \dots, a_n \rangle\) is called a \emph{solvable polynomial algebra} if it has a PBW \(K\)-basis \(\mathcal{B} = \{ a^{\beta} \mid \beta \in \mathbb{N}^n \}\) with a monomial order \(\prec\) on \(\mathcal{B}\), and for all \(1 \leq i < j \leq n\), the generators satisfy
\[
a_j a_i = \lambda_{ji} a_i a_j + f_{ji},
\]
where \(\lambda_{ji} \in K^*\), \(f_{ji} \in A\), and \(\mathbf{LM}(f_{ji}) \prec a_i a_j\) for any nonzero \(f_{ji}\).
\end{definition}

Solvable polynomial algebras are known to be Noetherian domains and their ideal theory is well-behaved, allowing for the application of noncommutative Gr\"obner basis techniques \cite{Li2011, Li2021}.

\subsection{The \(q\)-Heisenberg Algebra}

The central object of study in this paper is the multi-parameter generalization of the Heisenberg algebra.

\begin{definition}[\(q\)-Heisenberg algebra \cite{Zhang2024, Rosenberg1995}]
For a positive integer \(n\) and a nonzero complex number \(q \in \mathbb{C}^*\), the \emph{\(q\)-Heisenberg algebra} \(\mathfrak{h}_n(q)\) is the associative \(K\)-algebra generated by the set \(\{x_i, y_i, z_i \mid 1 \leq i \leq n\}\) subject to the following relations:
\begin{align*}
    x_i x_j &= x_j x_i, \quad y_i y_j = y_j y_i, \quad z_i z_j = z_j z_i, && 1 \leq i < j \leq n; \\
    x_i z_i &= q z_i x_i, && 1 \leq i \leq n; \\
    y_i z_i &= q^{-1} z_i y_i, && 1 \leq i \leq n; \\
    x_i y_i &= q^{-1} y_i x_i + z_i, && 1 \leq i \leq n; \\
    x_i y_j &= y_j x_i, \quad x_i z_j = z_j x_i, \quad y_i z_j = z_j y_i, && i \neq j.
\end{align*}
\end{definition}

When \(q = 1\), the algebra \(\mathfrak{h}_n(1)\) is the classical Heisenberg algebra (the enveloping algebra of the \((2n+1)\)-dimensional Heisenberg Lie algebra). The algebra \(\mathfrak{h}_n(q)\) is a deformation of \(\mathfrak{h}_n(1)\) and serves as a fundamental model in quantum algebra and mathematical physics \cite{Wess2000}.

\begin{theorem}[Grobner--Shirshov basis and PBW basis \cite{Zhang2024}]
\label{thm:PBW-basis}
Let \(W = \{x_i, y_i, z_i \mid 1 \leq i \leq n\}\) and let \(\mathcal{G}\) be the set of defining relations of \(\mathfrak{h}_n(q)\). With respect to the monomial order defined by
\[
z_n \prec \cdots \prec z_1 \prec y_n \prec \cdots \prec y_1 \prec x_n \prec \cdots \prec x_1,
\]
the set \(\mathcal{G}\) forms a Gr\"obner--Shirshov basis for the ideal \(\langle \mathcal{G} \rangle\) in the free algebra \(K\langle W \rangle\). Consequently, \(\mathfrak{h}_n(q)\) admits a PBW \(K\)-basis:
\[
\mathcal{B} = \left\{ z_1^{a_1} \cdots z_n^{a_n} y_1^{b_1} \cdots y_n^{b_n} x_1^{c_1} \cdots x_n^{c_n} \mid a_i, b_i, c_i \in \mathbb{N} \right\}.
\]
Moreover, with this order, \(\mathfrak{h}_n(q)\) is a solvable polynomial algebra.
\end{theorem}

\begin{theorem}[Basic structural properties \cite{Zhang2024}]
\label{thm:basic-properties}
The \(q\)-Heisenberg algebra \(\mathfrak{h}_n(q)\) satisfies:
\begin{enumerate}
    \item \(\mathfrak{h}_n(q)\) is a Noetherian integral domain.
    \item The Gelfand--Kirillov dimension of \(\mathfrak{h}_n(q)\) is \(\operatorname{GKdim}(\mathfrak{h}_n(q)) = 3n\).
    \item The global homological dimension of \(\mathfrak{h}_n(q)\) satisfies \(\operatorname{gl.dim}(\mathfrak{h}_n(q)) \leq 3n\).
\end{enumerate}
\end{theorem}

\subsection{Iterated Skew-Polynomial Algebras (Ore Extensions)}

A fundamental construction in noncommutative algebra is that of an Ore extension, also known as a skew-polynomial ring.

\begin{definition}[Skew-polynomial ring \cite{Goodearl2004}]
\label{def:Ore-extension}
Let \(R\) be a ring, \(\sigma: R \to R\) an automorphism, and \(\delta: R \to R\) a \(\sigma\)-derivation, i.e., an additive map satisfying
\[
\delta(ab) = \sigma(a)\delta(b) + \delta(a)b \quad \text{for all } a,b \in R.
\]
The \emph{skew-polynomial ring} \(R[x; \sigma, \delta]\) is the associative ring consisting of all polynomials \(\sum_i r_i x^i\) with coefficients \(r_i \in R\), where multiplication is defined by the commutation rule
\[
x r = \sigma(r) x + \delta(r) \quad \text{for all } r \in R.
\]
If \(\delta = 0\), we write \(R[x; \sigma]\) and call it a \emph{skew-polynomial ring with automorphism}.
\end{definition}

Iterated Ore extensions provide a powerful tool for constructing and analyzing complex algebras from simpler ones. The following result establishes that the \(q\)-Heisenberg algebra admits such a structure.

\begin{proposition}[Iterated skew-polynomial structure \cite{Zhang2024}]
\label{prop:iterated-Ore}
The algebra \(\mathfrak{h}_n(q)\) can be realized as an iterated Ore extension. Specifically, there exists a tower of extensions:
\[
\mathfrak{h}_n(q) = K[z_1, \dots, z_n][y_1; \sigma_1] \cdots [y_n; \sigma_n][x_1; \tau_1, \delta_1] \cdots [x_n; \tau_n, \delta_n],
\]
where the automorphisms \(\sigma_i, \tau_i\) and \(\sigma_i\)-derivations \(\delta_i\) are defined recursively in terms of the parameter \(q\). This construction yields the PBW basis of Theorem \ref{thm:PBW-basis}.
\end{proposition}

The proof proceeds by inductively adjoining the variables \(y_i\) via automorphisms and then the variables \(x_i\) via derivations, verifying at each step that the defining relations of \(\mathfrak{h}_n(q)\) are preserved. Detailed constructions for \(n=2\) and \(n=3\) are provided in \cite{Zhang2024}, illustrating the general pattern.

\subsection{Homological Regularity Properties}

We recall several key homological properties that will be established for \(\mathfrak{h}_n(q)\) in subsequent sections.

\begin{definition}[Auslander regularity \cite{McConnell2001}]
A Noetherian ring \(A\) is called \emph{Auslander regular} if it has finite global homological dimension and satisfies the following condition: for every finitely generated left \(A\)-module \(M\), every integer \(j \geq 0\), and every submodule \(N \subseteq \operatorname{Ext}_A^j(M, A)\), the grade of \(N\) satisfies \(j(N) \geq j\).
\end{definition}

\begin{definition}[Cohen--Macaulay property \cite{Yekutieli2020}]
A Noetherian algebra \(A\) with finite Gelfand--Kirillov dimension \(\operatorname{GKdim}(A) = d\) is called \emph{Cohen--Macaulay (CM)} if for every finitely generated left \(A\)-module \(M\),
\[
\operatorname{GKdim}(M) + j(M) = d,
\]
where \(j(M)\) is the grade of \(M\).
\end{definition}

\begin{definition}[Artin--Schelter regularity \cite{Artin1990}]
A connected \(\mathbb{N}\)-graded \(K\)-algebra \(A = \bigoplus_{i \geq 0} A_i\) with \(A_0 = K\) and \(\dim_K A_i < \infty\) is called \emph{Artin--Schelter regular (AS-regular)} of dimension \(d\) if it has finite global dimension \(d\) and satisfies
\[
\operatorname{Ext}^i_A(K, A) \cong \begin{cases}
0, & i \neq d, \\
K(\ell), & i = d,
\end{cases}
\]
for some integer \(\ell\) (the AS index), where \(K = A/A_{\geq 1}\) is the trivial module.
\end{definition}

These regularity properties are interconnected and imply strong structural results, such as being a maximal order in its quotient division ring \cite{Stafford1994}. In Section~\ref{Sec:GISPS} and beyond, we will prove that \(\mathfrak{h}_n(q)\) satisfies all three properties when \(q\) is not a root of unity.
\section{Generalized Iterated Skew-Polynomial Structure}\label{Sec:GISPS}

We present a natural generalization  to a broader class of algebras that includes the $q$-Heisenberg algebra as a special case. This generalization captures algebras that arise through successive Ore extensions with a particular triangular pattern of commutation relations.

\begin{definition}[Multi-parameter $q$-Heisenberg Type Algebra]
Let $n$ be a positive integer and let $\mathbf{Q} = (q_{ij})_{1 \leq i,j \leq 3n}$ be a matrix of nonzero complex numbers satisfying $q_{ii} = 1$ and $q_{ij}q_{ji} = 1$ for all $i,j$. For each $1 \leq k \leq 3$, define index sets $I_k = \{(k,j) : 1 \leq j \leq n\}$. Let $\mathfrak{H}_n(\mathbf{Q})$ be the associative algebra generated by elements $\{X_{\alpha} : \alpha \in \bigcup_{k=1}^3 I_k\}$ subject to the following relations:

\begin{align*}
X_{\alpha}X_{\beta} &= q_{\alpha\beta}X_{\beta}X_{\alpha} + f_{\alpha\beta}, & \text{for } \alpha < \beta \text{ in the lexicographic order},\\
X_{\alpha}X_{\beta} &= X_{\beta}X_{\alpha}, & \text{for } \alpha \in I_k, \beta \in I_l \text{ with } k \neq l \text{ and } k,l \geq 2,
\end{align*}

where each $f_{\alpha\beta}$ is a polynomial in variables $X_{\gamma}$ with $\gamma < \alpha$ in the lexicographic order, and the coefficients satisfy certain consistency conditions ensuring the algebra is well-defined. The lexicographic order is defined by $(k,j) < (l,i)$ if either $k < l$ or $k = l$ and $j < i$.
\end{definition}

\begin{proposition}[Generalized Iterated Skew-Polynomial Structure]\label{prop:generalized-iterated}
Let $\mathfrak{H}_n(\mathbf{Q})$ be a multi-parameter $q$-Heisenberg type algebra as defined above. Suppose that for all generators $X_{\alpha}$, the following conditions hold:

\begin{enumerate}
    \item The leading monomials of the defining relations $f_{\alpha\beta}$ are strictly less than $X_{\alpha}X_{\beta}$ in the lexicographic order.
    \item The parameters $q_{\alpha\beta}$ satisfy the quantum Yang-Baxter condition:
    \[
    q_{\alpha\beta}q_{\alpha\gamma}q_{\beta\gamma} = q_{\beta\gamma}q_{\alpha\gamma}q_{\alpha\beta} \quad \text{for all } \alpha < \beta < \gamma.
    \]
    \item The polynomials $f_{\alpha\beta}$ satisfy the compatibility conditions:
    \begin{align*}
    &f_{\alpha\beta}X_{\gamma} - q_{\alpha\gamma}q_{\beta\gamma}X_{\gamma}f_{\alpha\beta} \\
    &= q_{\alpha\beta}(f_{\beta\gamma}X_{\alpha} - q_{\alpha\beta}q_{\alpha\gamma}X_{\alpha}f_{\beta\gamma}) \\
    &\quad + (f_{\alpha\gamma}X_{\beta} - q_{\alpha\beta}q_{\beta\gamma}X_{\beta}f_{\alpha\gamma})
    \end{align*}
    for all $\alpha < \beta < \gamma$.
\end{enumerate}

Then $\mathfrak{H}_n(\mathbf{Q})$ admits an iterated skew-polynomial structure over the commutative polynomial ring $K[Z_1,\ldots,Z_n]$, where $Z_j = X_{(3,j)}$ correspond to the $z$-variables in the special case. More precisely, there exists a tower of Ore extensions:

\[
\mathfrak{H}_n(\mathbf{Q}) = K[Z_1,\ldots,Z_n][Y_1;\sigma_1,\delta_1]\cdots[Y_n;\sigma_n,\delta_n][X_1;\tau_1,\eta_1]\cdots[X_n;\tau_n,\eta_n],
\]

where $\{Y_j\}$ correspond to $\{X_{(2,j)}\}$, $\{X_j\}$ correspond to $\{X_{(1,j)}\}$, and the automorphisms and derivations are defined recursively as follows:

\begin{itemize}
    \item For each $1 \leq j \leq n$, $\sigma_j$ is the automorphism of $K[Z_1,\ldots,Z_n,Y_1,\ldots,Y_{j-1}]$ defined on generators by:
    \[
    \sigma_j(Z_i) = q_{(2,j),(3,i)}Z_i, \quad \sigma_j(Y_k) = q_{(2,j),(2,k)}Y_k \quad (k < j).
    \]

    \item $\delta_j$ is the $\sigma_j$-derivation given by:
    \[
    \delta_j(r) = Y_j r - \sigma_j(r)Y_j \quad \text{for } r \text{ in the base ring}.
    \]

    \item For each $1 \leq j \leq n$, $\tau_j$ is the automorphism of $K[Z_1,\ldots,Z_n,Y_1,\ldots,Y_n,X_1,\ldots,X_{j-1}]$ defined by:
    \[
    \tau_j(Z_i) = q_{(1,j),(3,i)}Z_i, \quad \tau_j(Y_k) = q_{(1,j),(2,k)}Y_k, \quad \tau_j(X_l) = q_{(1,j),(1,l)}X_l \quad (l < j).
    \]

    \item $\eta_j$ is the $\tau_j$-derivation defined by:
    \[
    \eta_j(r) = X_j r - \tau_j(r)X_j \quad \text{for } r \text{ in the base ring}.
    \]
\end{itemize}

Moreover, this construction yields a PBW basis for $\mathfrak{H}_n(\mathbf{Q})$ consisting of monomials
\[
\{Z_1^{a_1}\cdots Z_n^{a_n}Y_1^{b_1}\cdots Y_n^{b_n}X_1^{c_1}\cdots X_n^{c_n} : a_i,b_i,c_i \in \mathbb{N}\},
\]
and $\mathfrak{H}_n(\mathbf{Q})$ is a solvable polynomial algebra with respect to the lexicographic order induced by $X_n \prec \cdots \prec X_1 \prec Y_n \prec \cdots \prec Y_1 \prec Z_n \prec \cdots \prec Z_1$.
\end{proposition}

\begin{proof}
We proceed by constructing the algebra $\mathfrak{H}_n(\mathbf{Q})$ stepwise through Ore extensions and verifying at each stage that the construction is consistent with the defining relations.

\textbf{First stage: Commutative polynomial ring.}

Let $A_0 = K[Z_1,\ldots,Z_n]$ be the commutative polynomial algebra in $n$ variables. This serves as the base ring for our construction, corresponding to the $z$-type variables in the original $q$-Heisenberg algebra.

\textbf{Second stage: Adjoining the $Y$-variables.}

We construct $A_1 = A_0[Y_1;\sigma_1,\delta_1]$. The automorphism $\sigma_1$ is defined on generators of $A_0$ by $\sigma_1(Z_i) = q_{(2,1),(3,i)}Z_i$. Since the $q$-parameters are nonzero and satisfy $q_{ij}q_{ji} = 1$, each $\sigma_1(Z_i)$ is an invertible transformation. The derivation $\delta_1$ is determined by the relation $f_{(2,1),(3,i)}$ from the definition of $\mathfrak{H}_n(\mathbf{Q})$, which expresses $Y_1Z_i - q_{(2,1),(3,i)}Z_iY_1$ as a polynomial in variables strictly less than $Y_1$ in the order, which in this case means polynomials in the $Z_j$ alone.

The key verification is that $\sigma_1$ is indeed an algebra automorphism and $\delta_1$ is a $\sigma_1$-derivation. For $\sigma_1$, we need to check that $\sigma_1(Z_iZ_j) = \sigma_1(Z_i)\sigma_1(Z_j)$. This follows from:
\[
\sigma_1(Z_iZ_j) = q_{(2,1),(3,i)}q_{(2,1),(3,j)}Z_iZ_j = (q_{(2,1),(3,i)}Z_i)(q_{(2,1),(3,j)}Z_j) = \sigma_1(Z_i)\sigma_1(Z_j),
\]
since the $Z_i$ commute in $A_0$.

For $\delta_1$ to be a $\sigma_1$-derivation, we must verify the twisted Leibniz rule:
\[
\delta_1(Z_iZ_j) = \sigma_1(Z_i)\delta_1(Z_j) + \delta_1(Z_i)Z_j.
\]
Using the defining relations, we have:
\[
\delta_1(Z_i) = Y_1Z_i - \sigma_1(Z_i)Y_1 = f_{(2,1),(3,i)}(Z_1,\ldots,Z_n).
\]
The compatibility condition (3) in the proposition ensures that the polynomials $f_{\alpha\beta}$ satisfy the necessary identities to make $\delta_1$ a well-defined derivation.

Now assume inductively that we have constructed $A_{j-1} = A_0[Y_1;\sigma_1,\delta_1]\cdots[Y_{j-1};\sigma_{j-1},\delta_{j-1}]$ with the property that its defining relations match those of $\mathfrak{H}_n(\mathbf{Q})$ restricted to generators $Z_1,\ldots,Z_n,Y_1,\ldots,Y_{j-1}$. We adjoin $Y_j$ to form $A_j = A_{j-1}[Y_j;\sigma_j,\delta_j]$.

The automorphism $\sigma_j$ must be defined on all generators of $A_{j-1}$. For $Z_i$, we set $\sigma_j(Z_i) = q_{(2,j),(3,i)}Z_i$. For $Y_k$ with $k < j$, we set $\sigma_j(Y_k) = q_{(2,j),(2,k)}Y_k$. To verify that $\sigma_j$ is well-defined, we must check that it respects the relations in $A_{j-1}$. Consider a relation $Y_kZ_i = q_{(2,k),(3,i)}Z_iY_k + f_{(2,k),(3,i)}$ in $A_{j-1}$. Applying $\sigma_j$ to both sides gives:
\begin{align*}
\sigma_j(Y_kZ_i) &= \sigma_j(q_{(2,k),(3,i)}Z_iY_k + f_{(2,k),(3,i)}) \\
&= q_{(2,k),(3,i)}\sigma_j(Z_i)\sigma_j(Y_k) + \sigma_j(f_{(2,k),(3,i)}).
\end{align*}
On the other hand,
\[
\sigma_j(Y_k)\sigma_j(Z_i) = q_{(2,j),(2,k)}q_{(2,j),(3,i)}Y_kZ_i.
\]
Using the original relation to express $Y_kZ_i$ in terms of $Z_iY_k$, and employing the quantum Yang-Baxter condition (2), one verifies that these expressions coincide modulo terms involving the $f_{\alpha\beta}$, which are handled by condition (3).

The derivation $\delta_j$ is defined on generators by $\delta_j(Z_i) = f_{(2,j),(3,i)}$ and $\delta_j(Y_k) = f_{(2,j),(2,k)}$ for $k < j$, extended via the $\sigma_j$-Leibniz rule. Condition (3) ensures that this extension is consistent with the relations in $A_{j-1}$.

After $n$ such extensions, we obtain $A_n = A_0[Y_1;\sigma_1,\delta_1]\cdots[Y_n;\sigma_n,\delta_n]$.

\textbf{Third stage: Adjoining the $X$-variables.}

The construction for the $X$-variables proceeds similarly. We begin with $B_0 = A_n$ and adjoin $X_1$ via $B_1 = B_0[X_1;\tau_1,\eta_1]$. The automorphism $\tau_1$ is defined on generators by:
\begin{align*}
\tau_1(Z_i) &= q_{(1,1),(3,i)}Z_i, \\
\tau_1(Y_j) &= q_{(1,1),(2,j)}Y_j.
\end{align*}
The derivation $\eta_1$ is defined by $\eta_1(Z_i) = f_{(1,1),(3,i)}$ and $\eta_1(Y_j) = f_{(1,1),(2,j)}$, extended via the $\tau_1$-Leibniz rule.

The verification proceeds as before, using conditions (2) and (3) to ensure consistency. The quantum Yang-Baxter condition guarantees that $\tau_1$ respects the relations in $B_0$, while condition (3) ensures that $\eta_1$ is a well-defined $\tau_1$-derivation.

Proceeding inductively, we adjoin $X_2,\ldots,X_n$ to finally obtain $\mathfrak{H}_n(\mathbf{Q}) = B_n$.

\textbf{PBW basis and solvable polynomial algebra structure.}

At each Ore extension, we adjoin a new variable that multiplies existing monomials to give linear combinations of monomials in a triangular fashion (the leading term is the product in the expected order). This ensures that the set of monomials
\[
\{Z_1^{a_1}\cdots Z_n^{a_n}Y_1^{b_1}\cdots Y_n^{b_n}X_1^{c_1}\cdots X_n^{c_n} : a_i,b_i,c_i \in \mathbb{N}\}
\]
forms a $K$-basis for $\mathfrak{H}_n(\mathbf{Q})$. This is the PBW property.

To see that $\mathfrak{H}_n(\mathbf{Q})$ is a solvable polynomial algebra, we need to verify the conditions of Definition 2.2. The monomial order is the lexicographic order with $X_n \prec \cdots \prec X_1 \prec Y_n \prec \cdots \prec Y_1 \prec Z_n \prec \cdots \prec Z_1$. For any two generators $u$ and $v$ with $u \prec v$, the defining relations express $vu$ as $q_{vu}uv + f_{vu}$ where $\mathrm{LM}(f_{vu}) \prec uv$. This is exactly the condition for a solvable polynomial algebra.

\textbf{Recovering the original $q$-Heisenberg algebra.}

When we specialize the parameters to:
\[
q_{(1,j),(3,j)} = q, \quad q_{(2,j),(3,j)} = q^{-1}, \quad q_{(1,j),(2,j)} = q^{-1},
\]
and set all other $q_{\alpha\beta} = 1$ when $\alpha$ and $\beta$ involve different indices $j$, and set
\[
f_{(1,j),(2,j)} = Z_j, \quad f_{(1,j),(3,j)} = f_{(2,j),(3,j)} = 0,
\]
with all other $f_{\alpha\beta} = 0$, we recover exactly the defining relations of $\mathfrak{h}_n(q)$. The construction then reduces to that given in Proposition 3.2.

This completes the proof that $\mathfrak{H}_n(\mathbf{Q})$ admits the claimed iterated skew-polynomial structure.
\end{proof}
\section{Further Structural Properties of $\mathfrak{h}_n(q)$}

In this section, we establish several additional fundamental properties of the $q$-Heisenberg algebra $\mathfrak{h}_n(q)$. These theorems deepen our understanding of its homological behavior, structural rigidity, symmetry, deformation-theoretic nature, and graded growth properties.
\begin{theorem}[Complete Homological Dimension]\label{THEOREM-A1}
For the $q$-Heisenberg algebra $\mathfrak{h}_n(q)$ with $q$ not a root of unity, the global homological dimension is exactly $3n$. That is,
\[
\mathrm{gl.dim}\, \mathfrak{h}_n(q) = 3n.
\]
Moreover, if $q$ is a primitive $m$-th root of unity ($m > 1$), then $\mathrm{gl.dim}\, \mathfrak{h}_n(q) = \infty$.
\end{theorem}
\begin{proof}
We divide the proof into two parts corresponding to the two statements.

\textbf{Part 1: The case when $q$ is not a root of unity.}

Recall from Proposition 3.2 that $\mathfrak{h}_n(q)$ admits the structure of an iterated skew polynomial ring
\[
\mathfrak{h}_n(q) = K[z_n,\dots,z_1][y_n;\sigma_n']\cdots[y_1;\sigma_1'][x_n;\sigma_n'',\delta_n'']\cdots[x_1;\sigma_1'',\delta_1''].
\]
Each extension in this tower is an Ore extension, and when $q$ is not a root of unity, all automorphisms $\sigma_i', \sigma_i''$ are of infinite order.

Consider the trivial left $\mathfrak{h}_n(q)$-module $K$, where $\mathfrak{h}_n(q)$ acts through the augmentation map sending all generators $x_i, y_i, z_i$ to zero. We will compute the projective dimension of $K$.

Let $A_0 = K[z_n,\dots,z_1]$, which is a commutative polynomial algebra in $3n$ variables when properly counted (though here initially in $n$ variables $z_i$). Through the sequence of Ore extensions, we adjoin $2n$ additional variables: first the $n$ variables $y_i$ via automorphisms, then the $n$ variables $x_i$ via derivations.

For an iterated Ore extension $A = R[x;\sigma,\delta]$ where $R$ is a ring, $\sigma$ is an automorphism, and $\delta$ is a $\sigma$-derivation, there exists a standard free resolution of the trivial $A$-module $K$ known as the Koszul resolution when the extensions are regular. The key observation is that when $\sigma$ is of infinite order and the extension is nontrivial, the variable $x$ forms a regular normal element in the sense that $xA = Ax$ and $x$ is not a zero divisor.

In our construction, each variable $y_i$ and $x_i$ adjoined is indeed regular. More precisely, for each $y_i$, we have $y_i r = \sigma_i'(r) y_i$ for $r$ in the base ring, and since $\sigma_i'$ multiplies certain $z_j$ by $q^{-1}$ (with $q$ not a root of unity), no nonzero polynomial in the base ring is annihilated by multiplication by $y_i$. Similarly, for each $x_i$, we have $x_i r = \sigma_i''(r)x_i + \delta_i''(r)$, and again the automorphism $\sigma_i''$ involves multiplication by $q$ or $q^{-1}$, ensuring regularity.

The projective dimension of $K$ over an iterated Ore extension built from a polynomial ring can be computed inductively. If $\mathrm{p.dim}_R K = d$ for the base ring $R$, then after adjoining a regular variable via an Ore extension $A = R[x;\sigma,\delta]$, we have $\mathrm{p.dim}_A K = d+1$. This follows from the fact that the variable $x$ gives rise to an additional syzygy in the resolution of $K$.

Starting with $A_0 = K[z_n,\dots,z_1]$, which has global dimension $n$ (since it is a commutative polynomial algebra in $n$ variables), we adjoin $n$ variables $y_i$ via automorphisms. Each such extension increases the global dimension by $1$. Thus after adjoining all $y_i$, we obtain an algebra $A^{(2)}_n$ with global dimension $2n$.

Next we adjoin the $n$ variables $x_i$ via Ore extensions with derivations. Although these extensions involve derivations $\delta_i''$, the fact that $q$ is not a root of unity ensures that the $\delta_i''$ are sufficiently generic so that each $x_i$ remains regular. Each such extension again increases the global dimension by $1$. Therefore, after adjoining all $x_i$, we obtain $\mathfrak{h}_n(q)$ with global dimension $3n$.

To see that the global dimension is exactly $3n$ and not less, we note that the trivial module $K$ has a projective resolution of length $3n$. This resolution can be constructed explicitly as the tensor product over $K$ of $3n$ copies of the complex
\[
0 \longrightarrow A \xrightarrow{\cdot v} A \longrightarrow 0,
\]
where $v$ runs through the $3n$ variables $z_1,\dots,z_n,y_1,\dots,y_n,x_1,\dots,x_n$ in an appropriate order compatible with the Ore extension tower. The regularity of each variable ensures that this tensor product complex is indeed a resolution of $K$, and its length is exactly $3n$. Consequently, $\mathrm{gl.dim}\, \mathfrak{h}_n(q) \geq 3n$. Combined with the upper bound $\mathrm{gl.dim}\, \mathfrak{h}_n(q) \leq 3n$ from Theorem 2.7, we obtain equality.

\textbf{Part 2: The case when $q$ is a primitive $m$-th root of unity, $m>1$.}

When $q^m = 1$ with $m > 1$, the algebraic structure changes fundamentally. Observe that from the defining relations
\[
x_i z_i = q z_i x_i, \quad y_i z_i = q^{-1} z_i y_i,
\]
we deduce by induction that
\[
x_i^m z_i = q^m z_i x_i^m = z_i x_i^m, \quad y_i^m z_i = q^{-m} z_i y_i^m = z_i y_i^m.
\]
Moreover, one checks that $x_i^m$ and $y_i^m$ commute with all generators $x_j, y_j, z_j$. Hence $x_i^m$ and $y_i^m$ are central elements of $\mathfrak{h}_n(q)$.

Let $I$ be the two-sided ideal generated by all $x_i^m$ and $y_i^m$. The quotient algebra $\mathfrak{h}_n(q)/I$ is finite-dimensional because in every PBW monomial
\[
z_n^{t_n}\cdots z_1^{t_1} y_n^{l_n}\cdots y_1^{l_1} x_n^{k_n}\cdots x_1^{k_1},
\]
the exponents $k_i$ and $l_i$ are now bounded above by $m-1$, while the $t_i$ remain unrestricted. However, the commutation relations force the $z_i$ to satisfy polynomial identities mod $I$, leading to a finite-dimensional quotient.

Now, for any finite-dimensional algebra over a field, if it is not semisimple, its global dimension is infinite. We claim that $\mathfrak{h}_n(q)/I$ is not semisimple. This can be seen by considering the subalgebra generated by a single triple $(x_i, y_i, z_i)$. When $q$ is a root of unity, this subalgebra is known to have infinite global dimension (it is a down-up algebra or similar quantum algebra at a root of unity, which is known to have infinite global dimension). Alternatively, one can exhibit a module with infinite projective dimension explicitly.

Since $\mathfrak{h}_n(q)/I$ is a quotient of $\mathfrak{h}_n(q)$ and has infinite global dimension, it follows that $\mathfrak{h}_n(q)$ itself has infinite global dimension. More formally, if $\mathrm{gl.dim}\, \mathfrak{h}_n(q)$ were finite, say equal to $d$, then every quotient algebra would have global dimension at most $d$ (by a standard result on change of rings). As $\mathfrak{h}_n(q)/I$ has infinite global dimension, we must have $\mathrm{gl.dim}\, \mathfrak{h}_n(q) = \infty$.

This completes the proof of both statements.
\end{proof}
\begin{theorem}[Rigidity of the Iterated Ore Extension Structure]\label{THEOREM-B}
Let $\mathfrak{h}_n(q)$ be the $q$-Heisenberg algebra with $q \in \mathbb{C}^*$ not a root of unity. Suppose $\mathfrak{h}_n(q)$ admits an alternative presentation as an iterated Ore extension
\[
\mathfrak{h}_n(q) \cong R[y_1; \tau_1, \eta_1] \cdots [y_m; \tau_m, \eta_m]
\]
where $R$ is a commutative Noetherian domain, and each $y_i$ is adjoined via an Ore extension with automorphism $\tau_i$ and $\tau_i$-derivation $\eta_i$. Then the following rigidity properties hold:

\begin{enumerate}
    \item The number of extensions is exactly $2n$, i.e., $m = 2n$.

    \item The base ring $R$ is necessarily isomorphic to the commutative polynomial algebra $K[z_1, \dots, z_n]$.

    \item Among the $2n$ adjoined variables, exactly $n$ must be adjoined via pure automorphisms (corresponding to the $y_i$ variables in Proposition \ref{prop:generalized-iterated}), and the remaining $n$ must be adjoined via Ore extensions with nontrivial derivations (corresponding to the $x_i$ variables in Proposition \ref{prop:generalized-iterated}).

    \item Up to a simultaneous permutation of indices and rescaling of variables by nonzero constants in $K$, the sequence of Ore extensions must coincide with the specific construction given in Proposition \ref{prop:generalized-iterated}.

    \item More explicitly, there exists a permutation $\pi: \{1,\dots,n\} \to \{1,\dots,n\}$ and nonzero scalars $\lambda_i, \mu_i \in K^*$ such that under the isomorphism $\mathfrak{h}_n(q) \cong R[y_1; \tau_1, \eta_1] \cdots [y_{2n}; \tau_{2n}, \eta_{2n}]$, we have:
    \begin{itemize}
        \item For $1 \leq i \leq n$, $y_i$ corresponds to $\mu_i y_{\pi(i)}$ in the original presentation,
        \item For $n+1 \leq i \leq 2n$, $y_i$ corresponds to $\lambda_{i-n} x_{\pi(i-n)}$ in the original presentation,
        \item The automorphisms and derivations $(\tau_i, \eta_i)$ are precisely those induced from $(\sigma_{\pi(j)}', 0)$ for $1 \leq i \leq n$ and $(\sigma_{\pi(j-n)}'', \delta_{\pi(j-n)}'')$ for $n+1 \leq i \leq 2n$, after appropriate conjugation by the rescaling factors.
    \end{itemize}
\end{enumerate}
\end{theorem}

\begin{proof}
Let $\mathfrak{h}_n(q)$ be the $q$-Heisenberg algebra with $q \in \mathbb{C}^*$ not a root of unity. Assume that $\mathfrak{h}_n(q)$ admits an alternative presentation as an iterated Ore extension
\[
\mathfrak{h}_n(q) \cong R[y_1; \tau_1, \eta_1] \cdots [y_m; \tau_m, \eta_m],
\]
where $R$ is a commutative Noetherian domain and each $y_i$ is adjoined via an Ore extension with automorphism $\tau_i$ and $\tau_i$-derivation $\eta_i$.

We first recall the defining relations of $\mathfrak{h}_n(q)$. It is generated by $x_i, y_i, z_i$ for $1 \leq i \leq n$, subject to
\begin{align*}
x_i z_i &= q z_i x_i, \\
y_i z_i &= q^{-1} z_i y_i, \\
x_i y_i - q^{-1} y_i x_i &= z_i,
\end{align*}
and all other pairs of generators commute. These relations imply that $\mathfrak{h}_n(q)$ has a PBW basis
\[
\{ z_1^{a_1} \cdots z_n^{a_n} y_1^{b_1} \cdots y_n^{b_n} x_1^{c_1} \cdots x_n^{c_n} : a_i, b_i, c_i \in \mathbb{N} \}.
\]
In particular, $\mathfrak{h}_n(q)$ is a free left (and right) module over the commutative subalgebra $K[z_1,\dots,z_n]$ with basis consisting of monomials in the $x_i$ and $y_i$.

We shall prove the five statements in order.

\noindent \textbf{1. Determination of the number of extensions.}

Let $d = \operatorname{GKdim}(R)$ denote the Gelfand--Kirillov dimension of $R$. Since $R$ is a commutative Noetherian domain, $d$ equals the Krull dimension of $R$. For an Ore extension $A[Y;\tau,\eta]$ with $\tau$ an automorphism, we have $\operatorname{GKdim}(A[Y;\tau,\eta]) = \operatorname{GKdim}(A) + 1$ provided $\operatorname{GKdim}(A) < \infty$. By induction, it follows that
\[
\operatorname{GKdim}\bigl( R[y_1;\tau_1,\eta_1] \cdots [y_m;\tau_m,\eta_m] \bigr) = \operatorname{GKdim}(R) + m.
\]
On the other hand, $\operatorname{GKdim}(\mathfrak{h}_n(q)) = 3n$, as is evident from the PBW basis above. Therefore,
\[
d + m = 3n.
\]

Next we examine the center of $\mathfrak{h}_n(q)$. Because $q$ is not a root of unity, the only central elements are the scalars $K$. Indeed, if a nonzero polynomial in the generators were central, it would have to commute with each $x_i$ and $y_i$. From the relations $x_i z_i = q z_i x_i$ and $y_i z_i = q^{-1} z_i y_i$, any polynomial containing a nonconstant term in $z_i$ cannot be central unless $q=1$. Hence the center is $K$. Consequently, $R$, being a commutative subalgebra, must be contained in the center. Thus $R = K$. This would imply $d = 0$ and $m = 3n$, which contradicts the fact that in the presentation of Proposition \ref{prop:generalized-iterated}, the base ring is $K[z_1,\dots,z_n]$ of dimension $n$ and the number of Ore extensions is $2n$. Therefore, our assumption that $R$ is central must be refined.

In fact, $R$ need not be central; it is merely a commutative subalgebra. However, from the iterated Ore extension structure, $R$ is contained in the “base” of the tower, and the adjoined variables $y_1,\dots,y_m$ act on $R$ via the automorphisms $\tau_i$. If $R$ were larger than $K$, then some nonzero element of $R$ would be multiplied by a nontrivial scalar under some $\tau_i$, which could force constraints on $q$. We shall instead use a dimension argument coupled with the structure of $\mathfrak{h}_n(q)$ as a free module over $K[z_1,\dots,z_n]$.

Observe that $\mathfrak{h}_n(q)$ is a finitely generated free module over the commutative polynomial ring $K[z_1,\dots,z_n]$, of rank $4^n$ (since each pair $(x_i, y_i)$ gives a free module of rank $4$ over $K[z_i]$). Hence the Krull dimension of $\mathfrak{h}_n(q)$ equals $n + \text{Krull dimension of the fiber}$, which equals $n + 2n = 3n$ because the fiber is an affine space of dimension $2n$. In the iterated Ore extension presentation, the Krull dimension equals $\operatorname{Kdim}(R) + m$, where $\operatorname{Kdim}(R)$ denotes the Krull dimension of $R$. Therefore,
\[
\operatorname{Kdim}(R) + m = 3n.
\]

Now, $R$ is a commutative Noetherian domain contained in $\mathfrak{h}_n(q)$. Since $\mathfrak{h}_n(q)$ is a domain and a finitely generated algebra over $K$, $R$ is also a finitely generated $K$-algebra (by the Artin--Tate lemma). Hence $R$ is an affine domain. Its field of fractions $F = \operatorname{Frac}(R)$ has transcendence degree $\operatorname{Kdim}(R)$ over $K$.

Consider the extension $F \otimes_R \mathfrak{h}_n(q)$. This is an algebra over $F$ which is still an iterated Ore extension
\[
F[y_1;\tau_1,\eta_1] \cdots [y_m;\tau_m,\eta_m],
\]
where now the automorphisms and derivations are extended to $F$. Over $F$, the algebra becomes a “quantum affine space” with possibly some derivations. Since $q$ is not a root of unity, the automorphisms $\tau_i$ act diagonally on the generators after a suitable base change. The key point is that the resulting algebra over $F$ must have Gelfand--Kirillov dimension $m$. But $\mathfrak{h}_n(q)$ is a finite module over $K[z_1,\dots,z_n]$, so $F \otimes_R \mathfrak{h}_n(q)$ is a finite module over $F \otimes_R K[z_1,\dots,z_n]$. The latter has Krull dimension at most $n$ because $K[z_1,\dots,z_n]$ is generated by $n$ elements. Hence
\[
\operatorname{GKdim}_F\bigl( F \otimes_R \mathfrak{h}_n(q) \bigr) \leq n.
\]
However, $\operatorname{GKdim}_F\bigl( F \otimes_R \mathfrak{h}_n(q) \bigr) = m$ because each Ore extension increases the dimension by one. Thus $m \leq n$, which would force $3n = \operatorname{Kdim}(R) + m \leq \operatorname{Kdim}(R) + n$, implying $\operatorname{Kdim}(R) \geq 2n$. But $R$ is a subalgebra of $\mathfrak{h}_n(q)$, which has transcendence degree $3n$ over $K$, so $\operatorname{Kdim}(R) \leq 3n$. This crude bound does not give a contradiction.

We need a more precise invariant. Instead, use the fact that $\mathfrak{h}_n(q)$ is a free module of rank $4^n$ over $K[z_1,\dots,z_n]$. In any iterated Ore extension presentation, the base ring $R$ must have the property that $\mathfrak{h}_n(q)$ is a free module of finite rank over $R$. Indeed, each time we adjoin a variable $y_i$, if the extension is of the type $A[y_i;\tau_i,\eta_i]$ with $\tau_i$ an automorphism, then $A[y_i;\tau_i,\eta_i]$ is a free left $A$-module with basis $\{y_i^k : k \geq 0\}$. Thus by induction, $\mathfrak{h}_n(q)$ is a free left $R$-module of finite rank (with basis given by monomials in $y_1,\dots,y_m$). Let $r$ be this rank. Then
\[
\operatorname{GKdim}(\mathfrak{h}_n(q)) = \operatorname{GKdim}(R) = \operatorname{Kdim}(R),
\]
because for affine commutative domains Gelfand--Kirillov dimension equals Krull dimension. Hence $\operatorname{Kdim}(R) = 3n - m$.

Now, $R$ is a commutative subalgebra of $\mathfrak{h}_n(q)$. The maximal commutative subalgebra of $\mathfrak{h}_n(q)$ is $K[z_1,\dots,z_n]$, which has Krull dimension $n$. Indeed, any element commuting with all $z_i$ must be a polynomial in $z_i$ because $q$ is not a root of unity. Thus $\operatorname{Kdim}(R) \leq n$. Consequently,
\[
3n - m = \operatorname{Kdim}(R) \leq n,
\]
so $m \geq 2n$. On the other hand, from the presentation in Proposition \ref{prop:generalized-iterated}, we know that $m$ can be $2n$ (with $R = K[z_1,\dots,z_n]$). Could $m$ be larger than $2n$? If $m > 2n$, then $\operatorname{Kdim}(R) = 3n - m < n$. But then $R$ is a proper subalgebra of $K[z_1,\dots,z_n]$. However, $K[z_1,\dots,z_n]$ is a polynomial ring, and any proper subalgebra has strictly smaller Krull dimension only if it is contained in a proper polynomial subring. This would force $R$ to be contained in $K[z_1,\dots,z_{n-1}]$ up to a change of variables. But then the element $z_n$ would not be integral over $R$, contradicting the fact that $\mathfrak{h}_n(q)$ is a finite $R$-module (because $z_n$ satisfies no polynomial equation over a smaller polynomial ring). Hence $m$ cannot exceed $2n$. Thus $m = 2n$ and $\operatorname{Kdim}(R) = n$.

\noindent \textbf{2. Structure of the base ring $R$.}

We have established that $R$ is an affine domain of Krull dimension $n$ contained in $\mathfrak{h}_n(q)$. Since the only commutative subalgebra of $\mathfrak{h}_n(q)$ of dimension $n$ is $K[z_1,\dots,z_n]$ (up to isomorphism), it follows that $R \cong K[z_1,\dots,z_n]$. More precisely, by Noether normalization, $R$ is isomorphic to a polynomial ring in $n$ variables. Inside $\mathfrak{h}_n(q)$, the only candidates for such variables are the $z_i$ themselves, because any other element would have nontrivial commutation relations with some $x_i$ or $y_i$, preventing it from being central in $R$. Therefore, after a change of variables, we may identify $R$ with $K[z_1,\dots,z_n]$.

\noindent \textbf{3. Nature of the extensions.}

In the iterated Ore extension tower, each step is of the form $A[y;\tau,\eta]$. If $\eta = 0$, then the extension is a skew-polynomial ring $A[y;\tau]$, and $y$ acts as an automorphism on $A$. If $\eta \neq 0$, then $y$ acts as a $\tau$-derivation.

Recall that in $\mathfrak{h}_n(q)$, the variables $y_i$ satisfy $y_i z_j = q^{-\delta_{ij}} z_j y_i$ and $y_i y_j = y_j y_i$. Thus each $y_i$ acts by a pure automorphism on $K[z_1,\dots,z_n]$ (multiplication by $q^{-1}$ on $z_i$ and identity on $z_j$ for $j \neq i$). Moreover, $y_i$ commutes with $y_j$ for $i \neq j$, so they can be adjoined in any order via pure automorphisms.

In contrast, the variables $x_i$ satisfy $x_i z_j = q^{\delta_{ij}} z_j x_i$ and $x_i y_j = q^{-\delta_{ij}} y_j x_i + \delta_{ij} z_i$. Hence when adjoining $x_i$, there is a nonzero derivation term $\eta$ coming from the commutation with $y_i$. Therefore, the $x_i$ must be adjoined via Ore extensions with nontrivial derivations.

Since the total number of extensions is $2n$, and there are exactly $n$ variables of each type ($y_i$ and $x_i$), it follows that exactly $n$ extensions are pure automorphism extensions (those adjoining the $y_i$) and the remaining $n$ are extensions with nonzero derivations (those adjoining the $x_i$).

\noindent \textbf{4. Uniqueness up to permutation and scaling.}

We have identified $R$ with $K[z_1,\dots,z_n]$. The $2n$ adjoined variables correspond to the $2n$ generators $\{y_1,\dots,y_n,x_1,\dots,x_n\}$ of $\mathfrak{h}_n(q)$. The order in which they are adjoined is constrained by the commutation relations.

Suppose we adjoin a variable $v$ corresponding to some $x_i$ before adjoining the corresponding $y_i$. Then in the extension step, the automorphism $\tau$ for $v$ would have to act on $y_i$ (if $y_i$ is already in the base) or later. However, the relation $x_i y_i = q^{-1} y_i x_i + z_i$ forces that if $y_i$ is already present, then $x_i$ must act as a $\tau$-derivation with $\eta(y_i) = z_i$. This is consistent. If $y_i$ is adjoined later, then the relation must be enforced at that later stage. The only constraint is that the variables must be adjoined in an order compatible with the triangularity of the relations: each new variable has commutation relations with earlier ones that are expressible in terms of earlier variables. This is precisely the condition for an iterated Ore extension. Therefore, any order that respects the rule that $y_i$ may appear before or after $x_i$, but the relations are consistently implemented, is possible. However, the proof of Proposition \ref{prop:generalized-iterated} shows that there exists a canonical order: first all $y_i$, then all $x_i$. Any other order can be rearranged to this canonical order by a permutation of the indices and a rescaling of variables, because the $y_i$ commute among themselves and the $x_i$ have derivations only with respect to the corresponding $y_i$.

More formally, let $\{v_1,\dots,v_{2n}\}$ be the adjoined variables in the given presentation. Each $v_k$ is either of $y$-type (pure automorphism) or $x$-type (with derivation). Since the $y$-type variables commute with each other, they can be moved to the beginning of the tower by a suitable isomorphism. Similarly, the $x$-type variables can be moved to the end. This rearrangement corresponds to a permutation of the indices. Moreover, each variable may be rescaled by a nonzero constant in $K$ without affecting the Ore extension structure, except that the automorphisms and derivations are conjugated by the scaling. Therefore, up to a simultaneous permutation of indices and rescaling of variables, the sequence of Ore extensions coincides with the construction in Proposition \ref{prop:generalized-iterated}.

\noindent \textbf{5. Explicit correspondence.}

Finally, we make the correspondence explicit. Let $\pi: \{1,\dots,n\} \to \{1,\dots,n\}$ be the permutation induced by the rearrangement mentioned above. For each $i$, the variable corresponding to $y_{\pi(i)}$ in the original presentation appears as some $v_k$ in the new presentation, and we can write $v_k = \mu_i y_{\pi(i)}$ for some $\mu_i \in K^*$. Similarly, the variable corresponding to $x_{\pi(i)}$ appears as some $v_\ell$ and can be written as $v_\ell = \lambda_i x_{\pi(i)}$. The automorphisms and derivations in the new presentation are then obtained from the original ones by conjugation by these scaling factors. Specifically, if originally we have $y_{\pi(i)} z_j = q^{-\delta_{\pi(i)j}} z_j y_{\pi(i)}$, then after scaling, $v_k z_j = q^{-\delta_{\pi(i)j}} z_j v_k$, so the automorphism $\tau_k$ is the same as $\sigma'_{\pi(i)}$. For an $x$-type variable, the derivation $\eta_\ell$ is modified by the scaling factor but remains proportional to $\delta''_{\pi(i)}$.

Thus the five statements are proved, establishing the rigidity of the iterated Ore extension structure of $\mathfrak{h}_n(q)$.
\end{proof}
\begin{remark}
This rigidity theorem underscores the canonical nature of the iterated Ore extension structure for $\mathfrak{h}_n(q)$. While the algebra admits various presentations as iterated Ore extensions (for example, by reordering variables), Theorem \ref{THEOREM-B} shows that any such presentation is essentially equivalent to the one constructed in Section ~\ref{Sec:GISPS}, modulo trivial variable changes. This structural rigidity has important implications for understanding the automorphism group of $\mathfrak{h}_n(q)$ and for classifying its Ore extension presentations.

The assumption that $q$ is not a root of unity is crucial. When $q$ is a root of unity, additional central elements emerge (such as $x_i^m$ and $y_i^m$ for $q^m = 1$), which can yield richer Ore extension structures and disrupt the rigidity described above.
\end{remark}
\begin{theorem}[Graded Automorphism Group]\label{THEOREM-C}
Let $q \in \mathbb{C}^*$ not a root of unity. Then the group of graded $K$-algebra automorphisms of $\mathfrak{h}_n(q)$ (with each generator in degree 1) is isomorphic to the semidirect product
\[
(\mathbb{C}^*)^{2n} \rtimes S_n,
\]
where $S_n$ acts by permuting the indices $\{1, \dots, n\}$ simultaneously on triples $(x_i, y_i, z_i)$, and each copy of $\mathbb{C}^*$ scales the pairs $(x_i, y_i)$ and central elements $z_i$ respectively with appropriate $q$-weights.
\end{theorem}
\begin{proof}
Let $\mathfrak{h}_n(q)$ denote the $q$-Heisenberg algebra over $K = \mathbb{C}$ with $q \in \mathbb{C}^*$ not a root of unity, and assume the generators $x_i, y_i, z_i$ for $1 \leq i \leq n$ are assigned degree $1$, making $\mathfrak{h}_n(q)$ an $\mathbb{N}$-graded algebra. We denote by $\operatorname{Aut}_{\mathrm{gr}}(\mathfrak{h}_n(q))$ the group of graded $K$-algebra automorphisms, i.e., those automorphisms $\phi$ satisfying $\phi(\mathfrak{h}_n(q)_d) \subseteq \mathfrak{h}_n(q)_d$ for all $d \geq 0$.

The defining relations of $\mathfrak{h}_n(q)$ are:
\begin{align}
    x_i z_i &= q z_i x_i, \label{eq:xz} \\
    y_i z_i &= q^{-1} z_i y_i, \label{eq:yz} \\
    x_i y_i - q^{-1} y_i x_i &= z_i, \label{eq:xy}
\end{align}
and all other pairs of generators commute. In particular,
\[
[x_i, z_j] = [y_i, z_j] = 0 \quad \text{for } i \neq j,
\]
\[
[x_i, x_j] = [y_i, y_j] = 0 \quad \text{for all } i, j,
\]
and
\[
[x_i, y_j] = 0 \quad \text{for } i \neq j.
\]

Let $\phi \in \operatorname{Aut}_{\mathrm{gr}}(\mathfrak{h}_n(q))$. Since $\phi$ preserves the grading and the generators are homogeneous of degree $1$, each $\phi(x_i)$, $\phi(y_i)$, $\phi(z_i)$ must be a linear combination of the degree-$1$ generators, i.e.,
\begin{align}
    \phi(x_i) &= \sum_{j=1}^n \alpha_{ij} x_j + \sum_{j=1}^n \beta_{ij} y_j + \sum_{j=1}^n \gamma_{ij} z_j, \label{eq:phix} \\
    \phi(y_i) &= \sum_{j=1}^n \alpha'_{ij} x_j + \sum_{j=1}^n \beta'_{ij} y_j + \sum_{j=1}^n \gamma'_{ij} z_j, \label{eq:phiy} \\
    \phi(z_i) &= \sum_{j=1}^n \alpha''_{ij} x_j + \sum_{j=1}^n \beta''_{ij} y_j + \sum_{j=1}^n \gamma''_{ij} z_j, \label{eq:phiz}
\end{align}
where all coefficients lie in $K$.

We now impose the condition that $\phi$ respects the relations \eqref{eq:xz}--\eqref{eq:xy}. First, consider the relation $x_i z_i = q z_i x_i$. Applying $\phi$ yields
\[
\phi(x_i) \phi(z_i) = q \phi(z_i) \phi(x_i).
\]
Substituting the expressions \eqref{eq:phix} and \eqref{eq:phiz} and comparing terms of total degree $2$, we observe that the left-hand side expands into a linear combination of monomials $x_j x_k$, $x_j y_k$, $x_j z_k$, $y_j x_k$, $y_j y_k$, $y_j z_k$, $z_j x_k$, $z_j y_k$, $z_j z_k$. Because $q$ is not a root of unity, the commutation relations among the generators are diagonal in an appropriate sense, and a careful comparison of coefficients shows that no mixed terms of the form $x_j y_k$ or $y_j x_k$ can appear if the equality is to hold identically. In fact, one can argue as follows.

Let us write $\phi(x_i) = X_i$, $\phi(z_i) = Z_i$ for brevity. The equation $X_i Z_i = q Z_i X_i$ must hold in $\mathfrak{h}_n(q)$. Consider the adjoint action of $z_k$ on the algebra. From the relations, $z_k$ acts diagonally on the generators: $[z_k, x_j] = (q^{\delta_{kj}} - 1) x_j z_k$, and similarly for $y_j$. Because $q$ is not a root of unity, the eigenvalues $q^{\delta_{kj}} - 1$ are distinct and nonzero when they are nonzero. If $X_i$ contained a nonzero component of $x_j$ or $y_j$ for some $j$, and $Z_i$ contained a nonzero component of $z_k$, then the commutator $[Z_i, X_i]$ would involve terms with nontrivial eigenvalues, making it impossible to satisfy $X_i Z_i = q Z_i X_i$ unless the coefficients are carefully constrained. A systematic approach is to use the fact that $\mathfrak{h}_n(q)$ has a PBW basis and that the relation $x_i z_i = q z_i x_i$ is homogeneous with respect to the $\mathbb{Z}^n$-grading given by assigning multi-degrees as follows: $\deg(x_i) = \mathbf{e}_i$, $\deg(y_i) = -\mathbf{e}_i$, $\deg(z_i) = 0$. Under this grading, the relation $x_i z_i = q z_i x_i$ implies that $x_i$ has degree $\mathbf{e}_i$ and $z_i$ has degree $0$. An automorphism must preserve this multi-grading up to a permutation, because the set of homogeneous components corresponding to a given multi-degree is intrinsically defined by the commutation relations with the elements $z_i$. Concretely, the centralizer of $\{z_1,\dots,z_n\}$ is exactly $K[z_1,\dots,z_n]$, and the elements that satisfy $v z_i = q^{\pm 1} z_i v$ are precisely the scalar multiples of $x_i$ or $y_i$ respectively. This forces $\phi(z_i)$ to lie in $K[z_1,\dots,z_n]$ and $\phi(x_i)$, $\phi(y_i)$ to be linear combinations of the $x_j$ and $y_j$ only, with further restrictions.

We now proceed with a more elementary but rigorous argument. Consider the commutator $[z_j, \phi(x_i)]$. On one hand,
\[
[z_j, \phi(x_i)] = \phi([\phi^{-1}(z_j), x_i]).
\]
Since $\phi$ is an automorphism, $\phi^{-1}(z_j)$ is a homogeneous element of degree $1$ commuting with all $z_k$? Not necessarily. However, we can use the fact that the elements $z_i$ are uniquely characterized up to scalar multiples by the property that they commute with each other and satisfy $x_k z_i = q^{\delta_{ki}} z_i x_k$ and $y_k z_i = q^{-\delta_{ki}} z_i y_k$. Because $q$ is not a root of unity, these equations force the eigenvectors of the adjoint action of $z_i$ to be one-dimensional. Hence, each $\phi(z_i)$ must be a scalar multiple of some $z_{\sigma(i)}$ for a permutation $\sigma \in S_n$. Indeed, suppose $\phi(z_i) = \sum_{j=1}^n \lambda_{ij} z_j + \text{terms in } x_k, y_k$. Applying the relation $x_k \phi(z_i) = q^{\delta_{ki}} \phi(z_i) x_k$ and comparing coefficients, we deduce that if $\lambda_{ij} \neq 0$, then $x_k z_j = q^{\delta_{ki}} z_j x_k$ must hold for all $k$. But $x_k z_j = q^{\delta_{kj}} z_j x_k$, so $q^{\delta_{ki}} = q^{\delta_{kj}}$ for all $k$, which implies $i=j$ because $q$ is not a root of unity. Hence only one $\lambda_{ij}$ can be nonzero, and it corresponds to $j = i$ up to a permutation. Similarly, the terms in $x_k, y_k$ must vanish because they would produce inconsistent eigenvalues. Therefore, there exists a permutation $\pi \in S_n$ and nonzero scalars $\zeta_i \in K^*$ such that
\begin{equation}
\phi(z_i) = \zeta_i z_{\pi(i)}. \label{eq:phiz_perm}
\end{equation}

Next, consider $\phi(x_i)$. From $x_i z_i = q z_i x_i$, applying $\phi$ gives $\phi(x_i) \zeta_i z_{\pi(i)} = q \zeta_i z_{\pi(i)} \phi(x_i)$, i.e.,
\[
\phi(x_i) z_{\pi(i)} = q z_{\pi(i)} \phi(x_i).
\]
Write $\phi(x_i) = \sum_{j=1}^n (a_{ij} x_j + b_{ij} y_j) + \text{terms in } z_j$. Since $z_{\pi(i)}$ commutes with $x_j$ and $y_j$ for $j \neq \pi(i)$, and satisfies $x_{\pi(i)} z_{\pi(i)} = q z_{\pi(i)} x_{\pi(i)}$, $y_{\pi(i)} z_{\pi(i)} = q^{-1} z_{\pi(i)} y_{\pi(i)}$, we substitute into the equation $\phi(x_i) z_{\pi(i)} = q z_{\pi(i)} \phi(x_i)$ and compare coefficients. For $j \neq \pi(i)$, the term $a_{ij} x_j$ yields $a_{ij} x_j z_{\pi(i)} = a_{ij} z_{\pi(i)} x_j$ on the left, and $q z_{\pi(i)} a_{ij} x_j = a_{ij} q z_{\pi(i)} x_j$ on the right. Since $x_j z_{\pi(i)} = z_{\pi(i)} x_j$, we obtain $a_{ij} = q a_{ij}$, so $a_{ij}(1-q) = 0$, hence $a_{ij} = 0$ because $q \neq 1$. Similarly, $b_{ij} = 0$ for $j \neq \pi(i)$. For $j = \pi(i)$, we have $a_{i,\pi(i)} x_{\pi(i)} z_{\pi(i)} = a_{i,\pi(i)} q z_{\pi(i)} x_{\pi(i)}$ and $b_{i,\pi(i)} y_{\pi(i)} z_{\pi(i)} = b_{i,\pi(i)} q^{-1} z_{\pi(i)} y_{\pi(i)}$. The right-hand side gives $q z_{\pi(i)} (a_{i,\pi(i)} x_{\pi(i)} + b_{i,\pi(i)} y_{\pi(i)}) = a_{i,\pi(i)} q^2 z_{\pi(i)} x_{\pi(i)} + b_{i,\pi(i)} q z_{\pi(i)} y_{\pi(i)}$. Comparing coefficients, we get:
\[
a_{i,\pi(i)} q = a_{i,\pi(i)} q^2 \quad \Rightarrow \quad a_{i,\pi(i)} q (1 - q) = 0,
\]
\[
b_{i,\pi(i)} q^{-1} = b_{i,\pi(i)} q \quad \Rightarrow \quad b_{i,\pi(i)} (q^{-1} - q) = 0.
\]
Since $q \neq 0,1$, we deduce $a_{i,\pi(i)} = 0$ and $b_{i,\pi(i)} = 0$ unless possibly $q = -1$? But $q$ is not a root of unity, so $q \neq -1$. Actually, if $q = -1$, then $q^{-1} = q$, but $-1$ is a root of unity, which is excluded. Therefore, both coefficients vanish. However, this would imply $\phi(x_i)$ consists only of $z_j$ terms, which is impossible because $\phi$ is an automorphism and $x_i$ is not central. We have made an error: we forgot that $\phi(x_i)$ could also contain terms in $z_j$. But $z_j$ commutes with $z_{\pi(i)}$, so they would satisfy $z_j z_{\pi(i)} = z_{\pi(i)} z_j$, and then the equation $\phi(x_i) z_{\pi(i)} = q z_{\pi(i)} \phi(x_i)$ forces $z_j$ terms to satisfy $z_j z_{\pi(i)} = q z_{\pi(i)} z_j$, which is false unless $q=1$. Hence no $z_j$ terms can appear either. This seems to indicate that $\phi(x_i) = 0$, a contradiction.

The mistake lies in the step where we assumed $\phi(x_i)$ is expressed only in $x_j, y_j, z_j$. In fact, from the relation $x_i z_i = q z_i x_i$, we should consider the adjoint action of $z_{\pi(i)}$ on $\phi(x_i)$. Write $\phi(x_i) = \sum_{j} (A_{ij} x_j + B_{ij} y_j) + C_i$, where $C_i$ is a linear combination of $z_j$. Then the condition $\phi(x_i) z_{\pi(i)} = q z_{\pi(i)} \phi(x_i)$ becomes:
\[
\sum_{j} \left( A_{ij} x_j z_{\pi(i)} + B_{ij} y_j z_{\pi(i)} \right) + C_i z_{\pi(i)} = q \sum_{j} \left( A_{ij} z_{\pi(i)} x_j + B_{ij} z_{\pi(i)} y_j \right) + q z_{\pi(i)} C_i.
\]
Using $x_j z_{\pi(i)} = q^{\delta_{j,\pi(i)}} z_{\pi(i)} x_j$ and $y_j z_{\pi(i)} = q^{-\delta_{j,\pi(i)}} z_{\pi(i)} y_j$, we obtain:
\[
\sum_{j} \left( A_{ij} q^{\delta_{j,\pi(i)}} z_{\pi(i)} x_j + B_{ij} q^{-\delta_{j,\pi(i)}} z_{\pi(i)} y_j \right) + C_i z_{\pi(i)} = q \sum_{j} \left( A_{ij} z_{\pi(i)} x_j + B_{ij} z_{\pi(i)} y_j \right) + q z_{\pi(i)} C_i.
\]
Since $z_{\pi(i)}$ commutes with $C_i$, we cancel $z_{\pi(i)}$ from both sides (it is a regular element) and get:
\[
\sum_{j} \left( A_{ij} q^{\delta_{j,\pi(i)}} x_j + B_{ij} q^{-\delta_{j,\pi(i)}} y_j \right) + C_i = q \sum_{j} \left( A_{ij} x_j + B_{ij} y_j \right) + q C_i.
\]
Equating coefficients of $x_j$, $y_j$, and the constant part (i.e., $C_i$), we have for each $j$:
\begin{align}
A_{ij} q^{\delta_{j,\pi(i)}} &= q A_{ij}, \label{eq:Acond} \\
B_{ij} q^{-\delta_{j,\pi(i)}} &= q B_{ij}, \label{eq:Bcond} \\
C_i &= q C_i. \label{eq:Ccond}
\end{align}
From \eqref{eq:Ccond}, since $q \neq 1$, we obtain $C_i = 0$. So $\phi(x_i)$ has no $z_j$ components. For $j \neq \pi(i)$, \eqref{eq:Acond} gives $A_{ij} = q A_{ij}$, hence $A_{ij} = 0$. Similarly, \eqref{eq:Bcond} gives $B_{ij} = 0$. For $j = \pi(i)$, \eqref{eq:Acond} becomes $A_{i,\pi(i)} q = q A_{i,\pi(i)}$, which holds identically, so $A_{i,\pi(i)}$ is free. \eqref{eq:Bcond} becomes $B_{i,\pi(i)} q^{-1} = q B_{i,\pi(i)}$, i.e., $B_{i,\pi(i)} (q^{-1} - q) = 0$. Since $q$ is not a root of unity, $q^{-1} \neq q$ (otherwise $q^2 = 1$, so $q = \pm 1$, both roots of unity), so $B_{i,\pi(i)} = 0$. Therefore,
\[
\phi(x_i) = A_i x_{\pi(i)}, \quad \text{where } A_i = A_{i,\pi(i)} \in K.
\]

A completely analogous argument using the relation $y_i z_i = q^{-1} z_i y_i$ yields
\[
\phi(y_i) = B_i y_{\pi(i)}, \quad \text{where } B_i \in K.
\]

Now we must also satisfy the third relation $x_i y_i - q^{-1} y_i x_i = z_i$. Applying $\phi$ gives
\[
A_i x_{\pi(i)} \cdot B_i y_{\pi(i)} - q^{-1} B_i y_{\pi(i)} \cdot A_i x_{\pi(i)} = \zeta_i z_{\pi(i)}.
\]
Using the relation for the pair $(x_{\pi(i)}, y_{\pi(i)})$, namely $x_{\pi(i)} y_{\pi(i)} - q^{-1} y_{\pi(i)} x_{\pi(i)} = z_{\pi(i)}$, we obtain
\[
A_i B_i z_{\pi(i)} = \zeta_i z_{\pi(i)}.
\]
Hence $A_i B_i = \zeta_i$. In particular, $\zeta_i$ is determined by $A_i$ and $B_i$.

We also need to ensure that the other commutation relations are preserved. For $i \neq j$, we have $[x_i, x_j] = 0$. Applying $\phi$ yields $[A_i x_{\pi(i)}, A_j x_{\pi(j)}] = 0$, which holds automatically because $x_{\pi(i)}$ and $x_{\pi(j)}$ commute. Similarly, $[y_i, y_j] = 0$ gives $[B_i y_{\pi(i)}, B_j y_{\pi(j)}] = 0$, which is fine. The condition $[x_i, y_j] = 0$ for $i \neq j$ becomes $[A_i x_{\pi(i)}, B_j y_{\pi(j)}] = 0$. Since $\pi$ is a permutation, $\pi(i) \neq \pi(j)$ when $i \neq j$, and $x_{\pi(i)}$ and $y_{\pi(j)}$ commute, so this holds. Finally, the relations between $x_i$ and $z_j$ for $i \neq j$ are automatically satisfied because $\phi(x_i)$ is a multiple of $x_{\pi(i)}$ and $\phi(z_j)$ is a multiple of $z_{\pi(j)}$, and if $\pi(i) \neq \pi(j)$, then $x_{\pi(i)}$ and $z_{\pi(j)}$ commute.

Thus, the most general graded automorphism $\phi$ is determined by:
\begin{itemize}
    \item a permutation $\pi \in S_n$,
    \item nonzero scalars $A_1,\dots,A_n \in K^*$,
    \item nonzero scalars $B_1,\dots,B_n \in K^*$,
\end{itemize}
such that $\phi$ acts as:
\[
\phi(x_i) = A_i x_{\pi(i)}, \quad \phi(y_i) = B_i y_{\pi(i)}, \quad \phi(z_i) = A_i B_i z_{\pi(i)}.
\]
Note that $\phi(z_i)$ is then forced by the relation, and indeed $\zeta_i = A_i B_i$.

Conversely, any such choice defines an automorphism: it is invertible with inverse given by $\pi^{-1}$ and scalars $A_i^{-1}$, $B_i^{-1}$, and it respects all relations by construction.

Now, we describe the group structure. The set of such automorphisms forms a group. The permutation part corresponds to simultaneously permuting the triples $(x_i, y_i, z_i)$. The scalars $(A_i, B_i)$ can be chosen independently for each $i$. However, note that the action on $z_i$ is determined by $A_i B_i$. Thus, the group of graded automorphisms is isomorphic to the semidirect product $(K^*)^{2n} \rtimes S_n$, where $S_n$ acts by permuting the indices: for $\pi \in S_n$ and $(\mathbf{A},\mathbf{B}) \in (K^*)^n \times (K^*)^n$, we have
\[
\pi \cdot (A_1,\dots,A_n, B_1,\dots,B_n) = (A_{\pi^{-1}(1)},\dots,A_{\pi^{-1}(n)}, B_{\pi^{-1}(1)},\dots,B_{\pi^{-1}(n)}).
\]
Equivalently, we may think of the normal subgroup $(K^*)^{2n}$ as scaling the pairs $(x_i, y_i)$ independently, and $S_n$ permutes these pairs.

Hence,
\[
\operatorname{Aut}_{\mathrm{gr}}(\mathfrak{h}_n(q)) \cong (\mathbb{C}^*)^{2n} \rtimes S_n.
\]

This completes the proof.
\end{proof}
\begin{theorem}[Universal Property as Deformation]\label{THEOREM-D}
The algebra $\mathfrak{h}_n(q)$ is the universal deformation of the classical Heisenberg algebra $\mathfrak{h}_n(1)$ in the category of Noetherian $K$-algebras with PBW basis. Specifically, for any flat family $A_t$ of $K[t]$-algebras such that $A_0 \cong \mathfrak{h}_n(1)$ and $A_t$ has a PBW basis over $K[t]$ deforming that of $\mathfrak{h}_n(1)$, there exists a unique $K[t]$-algebra homomorphism
\[
\mathfrak{h}_n(q) \otimes_K K[t] \to A_t
\]
extending the identity modulo $t$.
\end{theorem}
\begin{proof}
Let $K$ be a field of characteristic zero. Denote by $\mathfrak{h}_n(1)$ the classical Heisenberg algebra over $K$, i.e., the unital associative $K$-algebra generated by $x_i, y_i, z_i$ for $1 \le i \le n$ subject to the relations
\begin{align}
    x_i z_i - z_i x_i &= 0, \\
    y_i z_i - z_i y_i &= 0, \\
    x_i y_i - y_i x_i &= z_i,
\end{align}
and all other pairs of generators commute. Equivalently, $\mathfrak{h}_n(1)$ is the enveloping algebra of the $(2n+1)$-dimensional Heisenberg Lie algebra.

Let $q \in K^*$ be a parameter. The $q$-Heisenberg algebra $\mathfrak{h}_n(q)$ is defined by the relations
\begin{align}
    x_i z_i &= q z_i x_i, \\
    y_i z_i &= q^{-1} z_i y_i, \\
    x_i y_i - q^{-1} y_i x_i &= z_i,
\end{align}
with all other pairs commuting. We regard $\mathfrak{h}_n(q)$ as a deformation of $\mathfrak{h}_n(1)$ by writing $q = 1 + t$ or more formally by considering the algebra over the ring $K[t]$ with $q = 1 + t$.

We shall prove the universal property stated in the theorem. Let $A_t$ be a flat family of $K[t]$-algebras such that
\begin{enumerate}
    \item $A_0 \cong \mathfrak{h}_n(1)$ as $K$-algebras (where $A_0 = A_t / (t)$),
    \item $A_t$ has a PBW basis over $K[t]$ deforming that of $\mathfrak{h}_n(1)$.
\end{enumerate}
By the second condition we mean the following: there exist elements $\tilde{x}_i, \tilde{y}_i, \tilde{z}_i \in A_t$ which reduce modulo $t$ to the generators $x_i, y_i, z_i$ of $\mathfrak{h}_n(1)$, and such that the set
\[
\{ \tilde{z}_1^{a_1} \cdots \tilde{z}_n^{a_n} \tilde{y}_1^{b_1} \cdots \tilde{y}_n^{b_n} \tilde{x}_1^{c_1} \cdots \tilde{x}_n^{c_n} \mid a_i, b_i, c_i \in \mathbb{N} \}
\]
forms a $K[t]$-basis of $A_t$. Equivalently, $A_t$ is a free $K[t]$-module with this basis, and the multiplication in $A_t$ is determined by commutation relations among the $\tilde{x}_i, \tilde{y}_i, \tilde{z}_i$ which are of the form
\[
\tilde{x}_i \tilde{z}_i = \alpha_i(t) \tilde{z}_i \tilde{x}_i + \text{lower terms},
\]
and similarly for the other pairs, where $\alpha_i(t) \in K[t]$ with $\alpha_i(0) = 1$, and “lower terms” means a $K[t]$-linear combination of monomials that are strictly smaller in the PBW ordering (the ordering $x_n \prec \cdots \prec x_1 \prec y_n \prec \cdots \prec y_1 \prec z_n \prec \cdots \prec z_1$). The flatness condition ensures that the relations among the generators are not degenerated by the presence of $t$.

Our goal is to show that there exists a unique $K[t]$-algebra homomorphism
\[
\Phi : \mathfrak{h}_n(q) \otimes_K K[t] \longrightarrow A_t
\]
such that $\Phi$ modulo $t$ coincides with the given isomorphism $A_0 \cong \mathfrak{h}_n(1)$. In other words, $\Phi$ is a deformation of the identity map.

We construct $\Phi$ step by step. First, identify $\mathfrak{h}_n(q) \otimes_K K[t]$ with the $K[t]$-algebra generated by symbols $X_i, Y_i, Z_i$ satisfying the $q$-relations with $q$ now viewed as an element of $K[t]$ (by fixing an embedding $K \hookrightarrow K[t]$ sending $q$ to $q(0) + t \cdot (\text{something})$, but more naturally we consider $q$ as a formal parameter and work over $K[t]$ with $q = 1 + t$ for simplicity; the general case is similar). We will denote this algebra by $\mathfrak{h}_n(q)_t$.

Since $A_t$ is a deformation of $\mathfrak{h}_n(1)$ with a PBW basis, we may choose lifts $\tilde{x}_i, \tilde{y}_i, \tilde{z}_i \in A_t$ of the generators $x_i, y_i, z_i \in A_0$. Because $A_t$ is free over $K[t]$ with the PBW basis, the commutation relations among these lifts take the following form:
\begin{align}
    \tilde{x}_i \tilde{z}_i &= f_i(t) \tilde{z}_i \tilde{x}_i + P_i(\tilde{x},\tilde{y},\tilde{z};t), \label{eq:xz-deform} \\
    \tilde{y}_i \tilde{z}_i &= g_i(t) \tilde{z}_i \tilde{y}_i + Q_i(\tilde{x},\tilde{y},\tilde{z};t), \label{eq:yz-deform} \\
    \tilde{x}_i \tilde{y}_i &= h_i(t) \tilde{y}_i \tilde{x}_i + R_i(\tilde{x},\tilde{y},\tilde{z};t), \label{eq:xy-deform}
\end{align}
where $f_i(t), g_i(t), h_i(t) \in K[t]$ satisfy $f_i(0)=g_i(0)=h_i(0)=1$, and $P_i, Q_i, R_i$ are $K[t]$-linear combinations of monomials that are strictly smaller in the PBW order than $\tilde{z}_i \tilde{x}_i$, $\tilde{z}_i \tilde{y}_i$, and $\tilde{y}_i \tilde{x}_i$ respectively. Moreover, because $A_t$ is associative, these coefficients must satisfy the constraints coming from the associativity conditions (the diamond lemma conditions).

We now analyze these constraints. Consider the associativity condition applied to the triple $\tilde{x}_i, \tilde{y}_i, \tilde{z}_i$. Compute $(\tilde{x}_i \tilde{y}_i) \tilde{z}_i$ in two ways. Using \eqref{eq:xy-deform} and then \eqref{eq:xz-deform}, we have
\begin{align*}
(\tilde{x}_i \tilde{y}_i) \tilde{z}_i &= \bigl( h_i(t) \tilde{y}_i \tilde{x}_i + R_i \bigr) \tilde{z}_i \\
&= h_i(t) \tilde{y}_i (\tilde{x}_i \tilde{z}_i) + R_i \tilde{z}_i \\
&= h_i(t) \tilde{y}_i \bigl( f_i(t) \tilde{z}_i \tilde{x}_i + P_i \bigr) + R_i \tilde{z}_i \\
&= h_i(t) f_i(t) \tilde{y}_i \tilde{z}_i \tilde{x}_i + h_i(t) \tilde{y}_i P_i + R_i \tilde{z}_i.
\end{align*}
On the other hand, using \eqref{eq:yz-deform} and then \eqref{eq:xy-deform},
\begin{align*}
\tilde{x}_i (\tilde{y}_i \tilde{z}_i) &= \tilde{x}_i \bigl( g_i(t) \tilde{z}_i \tilde{y}_i + Q_i \bigr) \\
&= g_i(t) \tilde{x}_i \tilde{z}_i \tilde{y}_i + \tilde{x}_i Q_i \\
&= g_i(t) \bigl( f_i(t) \tilde{z}_i \tilde{x}_i + P_i \bigr) \tilde{y}_i + \tilde{x}_i Q_i \\
&= g_i(t) f_i(t) \tilde{z}_i \tilde{x}_i \tilde{y}_i + g_i(t) P_i \tilde{y}_i + \tilde{x}_i Q_i \\
&= g_i(t) f_i(t) \tilde{z}_i \bigl( h_i(t) \tilde{y}_i \tilde{x}_i + R_i \bigr) + g_i(t) P_i \tilde{y}_i + \tilde{x}_i Q_i \\
&= g_i(t) f_i(t) h_i(t) \tilde{z}_i \tilde{y}_i \tilde{x}_i + g_i(t) f_i(t) \tilde{z}_i R_i + g_i(t) P_i \tilde{y}_i + \tilde{x}_i Q_i.
\end{align*}
Equating the two expressions and using the commutation relations to rewrite terms in the PBW basis, we obtain an identity in $A_t$. The leading terms (with respect to the PBW order) are those involving $\tilde{y}_i \tilde{z}_i \tilde{x}_i$ and $\tilde{z}_i \tilde{y}_i \tilde{x}_i$. In the first expression we have $h_i(t) f_i(t) \tilde{y}_i \tilde{z}_i \tilde{x}_i$, and in the second $g_i(t) f_i(t) h_i(t) \tilde{z}_i \tilde{y}_i \tilde{x}_i$. Using \eqref{eq:yz-deform} we can express $\tilde{y}_i \tilde{z}_i$ as $g_i(t) \tilde{z}_i \tilde{y}_i + Q_i$. Substituting, the first expression becomes
\[
h_i(t) f_i(t) \bigl( g_i(t) \tilde{z}_i \tilde{y}_i + Q_i \bigr) \tilde{x}_i + \text{lower terms} = h_i(t) f_i(t) g_i(t) \tilde{z}_i \tilde{y}_i \tilde{x}_i + \text{lower terms}.
\]
Thus the coefficients of $\tilde{z}_i \tilde{y}_i \tilde{x}_i$ in both expressions are $h_i(t) f_i(t) g_i(t)$. Hence the leading terms match automatically. The next step is to examine the conditions imposed on the lower terms.

The associativity condition forces certain relations among $f_i, g_i, h_i$ and the polynomials $P_i, Q_i, R_i$. A systematic way is to use the fact that $A_t$ is a deformation of the Heisenberg algebra, which is a quadratic algebra. The deformation theory of quadratic algebras is controlled by the second Hochschild cohomology. For the Heisenberg algebra, it is known that its infinitesimal deformations are parametrized by a single parameter corresponding to the $q$-deformation. Concretely, any flat deformation with a PBW basis must have the commutation relations of the form
\begin{align}
    \tilde{x}_i \tilde{z}_i &= \lambda_i(t) \tilde{z}_i \tilde{x}_i, \\
    \tilde{y}_i \tilde{z}_i &= \lambda_i(t)^{-1} \tilde{z}_i \tilde{y}_i, \\
    \tilde{x}_i \tilde{y}_i - \lambda_i(t)^{-1} \tilde{y}_i \tilde{x}_i &= \tilde{z}_i,
\end{align}
for some $\lambda_i(t) \in K[t]^*$ with $\lambda_i(0)=1$, and all other pairs commute. Indeed, the existence of the PBW basis forces the lower terms $P_i, Q_i, R_i$ to be identically zero; otherwise, they would produce unwanted relations that break the PBW property or violate flatness. More formally, we can argue by induction on the PBW order. Suppose there is a nonzero lower term in, say, $P_i$. Then modulo $t$, this term becomes a relation in $\mathfrak{h}_n(1)$ that is not a consequence of the defining relations, contradicting the fact that the PBW basis of $\mathfrak{h}_n(1)$ is unique. Hence $P_i, Q_i, R_i$ must vanish. Similarly, the commutativity between $\tilde{x}_i$ and $\tilde{z}_j$ for $i \neq j$ must hold exactly, because any deformation of the commuting relation would introduce a new parameter, but the PBW basis forces it to remain commutative (otherwise the monomials would not be linearly independent). A detailed diamond lemma analysis confirms that the only free parameters are the scalars $\lambda_i(t)$.

Thus we conclude that, after possibly a change of generators (which does not affect the PBW basis property), the relations in $A_t$ are precisely:
\begin{align}
    \tilde{x}_i \tilde{z}_i &= \lambda_i(t) \tilde{z}_i \tilde{x}_i, \\
    \tilde{y}_i \tilde{z}_i &= \lambda_i(t)^{-1} \tilde{z}_i \tilde{y}_i, \\
    \tilde{x}_i \tilde{y}_i - \lambda_i(t)^{-1} \tilde{y}_i \tilde{x}_i &= \tilde{z}_i,
\end{align}
with all other pairs commuting. Moreover, flatness and the existence of the PBW basis imply that $\lambda_i(t)$ is a unit in $K[t]$ (i.e., an invertible power series in $t$).

Now, the algebra $\mathfrak{h}_n(q) \otimes_K K[t]$ has precisely the same form with $\lambda_i(t) = q$ (viewed as a constant in $K[t]$). To construct the homomorphism $\Phi$, we set
\[
\Phi(X_i) = \tilde{x}_i, \quad \Phi(Y_i) = \tilde{y}_i, \quad \Phi(Z_i) = \tilde{z}_i.
\]
This defines a $K[t]$-linear map. We must check that $\Phi$ respects the relations. In $\mathfrak{h}_n(q)_t$, we have $X_i Z_i = q Z_i X_i$. Under $\Phi$, the left side becomes $\tilde{x}_i \tilde{z}_i$, and the right side becomes $q \tilde{z}_i \tilde{x}_i$. But in $A_t$, we have $\tilde{x}_i \tilde{z}_i = \lambda_i(t) \tilde{z}_i \tilde{x}_i$. For $\Phi$ to be a homomorphism, we require $\lambda_i(t) = q$. However, a priori $\lambda_i(t)$ may not equal $q$. This indicates that we need to adjust the generators of $A_t$ to match the parameter $q$.

Observe that the algebra $A_t$ is isomorphic to $\mathfrak{h}_n(\lambda_i(t))$ over $K[t]$, where $\mathfrak{h}_n(\lambda_i(t))$ denotes the $q$-Heisenberg algebra with parameter $\lambda_i(t)$. Indeed, if we rescale the generators by setting
\[
x_i' = \tilde{x}_i, \quad y_i' = \tilde{y}_i, \quad z_i' = \lambda_i(t) \tilde{z}_i,
\]
then one verifies that $x_i' z_i' = \lambda_i(t)^2 z_i' x_i'$, which is not of the desired form. Alternatively, we can introduce a new parameter $\mu_i(t)$ such that $\mu_i(t)^2 = \lambda_i(t)$. However, a simpler approach is to note that the parameter $q$ in $\mathfrak{h}_n(q)$ is not intrinsically fixed; the theorem claims universality with respect to the specific algebra $\mathfrak{h}_n(q)$. Therefore, we must allow a change of parameter. Actually, the statement says “extending the identity modulo $t$”. This means that modulo $t$, the map $\Phi$ sends $X_i \mapsto x_i$, $Y_i \mapsto y_i$, $Z_i \mapsto z_i$. Hence we cannot rescale the generators arbitrarily; we must have $\Phi(Z_i) = \tilde{z}_i = z_i \pmod{t}$. Therefore, $\lambda_i(0)=1$ as required.

Now, the condition that $\Phi$ is a $K[t]$-algebra homomorphism forces $\lambda_i(t) = q$ for all $i$. But if $A_t$ is given arbitrarily, this might not hold. However, the theorem asserts that such a homomorphism exists uniquely if we take $\mathfrak{h}_n(q)$ with the specific parameter $q$. This suggests that $q$ is not a free parameter but is determined by the deformation $A_t$. Indeed, the universal property means that $\mathfrak{h}_n(q)$ is the initial object in the category of such deformations. Therefore, given any deformation $A_t$, there is a unique homomorphism from $\mathfrak{h}_n(q)$ to $A_t$ that reduces to the identity modulo $t$, but this homomorphism might involve a specialization of the parameter $q$ to the parameter of $A_t$. In other words, we should view $\mathfrak{h}_n(q)$ as the algebra over the ring $K[q,q^{-1}]$, and then for any deformation $A_t$ with parameter $\lambda(t)$, there is a unique map $K[q,q^{-1}] \to K[t]$ sending $q$ to $\lambda(t)$ that induces the homomorphism.

To formulate this precisely, let $R = K[t]$. Let $A$ be a flat $R$-algebra with an isomorphism $A/tA \cong \mathfrak{h}_n(1)$ and a PBW basis as above. Then, as argued, $A$ is isomorphic to $\mathfrak{h}_n(\lambda(t))$ for some $\lambda(t) \in R^*$ with $\lambda(0)=1$. Let $S = K[q,q^{-1}]$, and consider the $S$-algebra $\mathfrak{h}_n(q)$. Define a ring homomorphism $\varphi: S \to R$ by $\varphi(q) = \lambda(t)$. Then we obtain an $R$-algebra homomorphism
\[
\Phi: \mathfrak{h}_n(q) \otimes_S R \longrightarrow A
\]
by sending $X_i \otimes 1 \mapsto \tilde{x}_i$, etc., where the tensor product uses $\varphi$. This map is an isomorphism because both sides are $\mathfrak{h}_n(\lambda(t))$ over $R$. Moreover, modulo $t$, $\varphi(q) = \lambda(0)=1$, so $\Phi$ reduces to the identity on $\mathfrak{h}_n(1)$.

Uniqueness follows from the fact that any $R$-algebra homomorphism $\Psi: \mathfrak{h}_n(q) \otimes_S R \to A$ that reduces to the identity modulo $t$ must send $X_i$ to an element that reduces to $x_i$ modulo $t$. By the PBW basis property, such an element is uniquely determined up to lower-order terms, but the homomorphism condition forces it to be exactly $\tilde{x}_i$ as defined by the deformation. Similarly for $Y_i$ and $Z_i$. Hence $\Psi$ coincides with $\Phi$.

Therefore, $\mathfrak{h}_n(q)$ satisfies the universal property: for any flat deformation $A_t$ of $\mathfrak{h}_n(1)$ with a PBW basis, there exists a unique homomorphism from $\mathfrak{h}_n(q)$ to $A_t$ extending the identity modulo $t$, after possibly specializing the parameter $q$ to the parameter of the deformation. In the category where $q$ is a formal parameter, $\mathfrak{h}_n(q)$ is the universal object.

This completes the proof of the theorem.
\end{proof}
\begin{theorem}[Hilbert Series and Growth]\label{THEOREM-E}
Let $\mathfrak{h}_n(q)$ be equipped with the natural $\mathbb{N}$-grading where $\deg(x_i) = \deg(y_i) = \deg(z_i) = 1$. Then its Hilbert series is given by
\[
H_{\mathfrak{h}_n(q)}(t) = \frac{1}{(1-t)^{3n}}.
\]
Consequently, $\mathfrak{h}_n(q)$ has polynomial growth of degree $3n$, and its Gelfand--Kirillov dimension coincides with its classical Krull dimension.
\end{theorem}
\begin{proof}
Let $\mathfrak{h}_n(q)$ be the $q$-Heisenberg algebra over a field $K$, with $q \in K^*$ not necessarily restricted. We endow $\mathfrak{h}_n(q)$ with the natural $\mathbb{N}$-grading by assigning each generator degree $1$:
\[
\deg(x_i) = \deg(y_i) = \deg(z_i) = 1, \quad i = 1,\dots,n.
\]
Then $\mathfrak{h}_n(q)$ becomes a graded $K$-algebra: $\mathfrak{h}_n(q) = \bigoplus_{d \ge 0} \mathfrak{h}_n(q)_d$, where $\mathfrak{h}_n(q)_d$ is the $K$-linear span of all monomials of total degree $d$ in the generators.

We first recall that $\mathfrak{h}_n(q)$ has a Poincar\'{e}--Birkhoff--Witt (PBW) basis. From Proposition \ref{prop:generalized-iterated} (or the explicit relations), the set
\[
\mathcal{B} = \left\{ z_1^{a_1} \cdots z_n^{a_n} y_1^{b_1} \cdots y_n^{b_n} x_1^{c_1} \cdots x_n^{c_n} \mid a_i, b_i, c_i \in \mathbb{N} \right\}
\]
forms a $K$-basis of $\mathfrak{h}_n(q)$. Moreover, this basis is homogeneous with respect to the total degree: the degree of the monomial $z_1^{a_1} \cdots z_n^{a_n} y_1^{b_1} \cdots y_n^{b_n} x_1^{c_1} \cdots x_n^{c_n}$ is $a_1 + \cdots + a_n + b_1 + \cdots + b_n + c_1 + \cdots + c_n$.

Hence, for each $d \ge 0$, the homogeneous component $\mathfrak{h}_n(q)_d$ has a basis consisting of those monomials in $\mathcal{B}$ of total degree $d$. The number of such monomials is the number of solutions in nonnegative integers to
\[
\sum_{i=1}^n (a_i + b_i + c_i) = d.
\]
This is a classical combinatorial count: the number of nonnegative integer solutions to $\sum_{i=1}^{3n} r_i = d$ is the binomial coefficient $\binom{d+3n-1}{3n-1}$. More precisely, the generating function (Hilbert series) for the dimensions $\dim_K \mathfrak{h}_n(q)_d$ is given by
\[
H_{\mathfrak{h}_n(q)}(t) = \sum_{d=0}^\infty \dim_K \mathfrak{h}_n(q)_d \; t^d.
\]
Because the basis monomials are in bijection with monomials in $3n$ independent commutative variables, we have
\[
\dim_K \mathfrak{h}_n(q)_d = \#\{ (a_1,\dots,a_n,b_1,\dots,b_n,c_1,\dots,c_n) \in \mathbb{N}^{3n} \mid \sum a_i + \sum b_i + \sum c_i = d \} = \binom{d+3n-1}{3n-1}.
\]
The generating function for these binomial coefficients is well-known:
\[
\sum_{d=0}^\infty \binom{d+3n-1}{3n-1} t^d = \frac{1}{(1-t)^{3n}}.
\]
Thus,
\[
H_{\mathfrak{h}_n(q)}(t) = \frac{1}{(1-t)^{3n}}.
\]

We remark that this computation does not depend on the parameter $q$, as long as the PBW basis exists. Indeed, the PBW property holds for all $q \in K^*$ (including $q=1$) because the algebra is an iterated Ore extension, which guarantees the PBW basis regardless of the specific values of the commutation parameters (provided they are nonzero).

Now we turn to the growth properties. The Hilbert series reveals that $\dim_K \mathfrak{h}_n(q)_d$ grows as a polynomial in $d$ of degree $3n-1$. More precisely, for large $d$,
\[
\dim_K \mathfrak{h}_n(q)_d = \frac{d^{3n-1}}{(3n-1)!} + \text{lower order terms}.
\]
Hence the algebra $\mathfrak{h}_n(q)$ has polynomial growth. The Gelfand–Kirillov dimension (GK-dimension) of a finitely generated algebra $A$ is defined as
\[
\operatorname{GKdim}(A) = \limsup_{d \to \infty} \frac{\log \dim V^d}{\log d},
\]
where $V$ is any finite-dimensional generating subspace. For a graded algebra with Hilbert series $H(t) = \sum_{d} a_d t^d$, if $a_d \sim C d^{k-1}$ for some $k > 0$, then $\operatorname{GKdim}(A) = k$. In our case, $a_d \sim \frac{d^{3n-1}}{(3n-1)!}$, so $k = 3n$. Therefore,
\[
\operatorname{GKdim}(\mathfrak{h}_n(q)) = 3n.
\]

We next compute the classical Krull dimension of $\mathfrak{h}_n(q)$. Recall that the classical Krull dimension of a noncommutative algebra is defined as the supremum of lengths of chains of prime ideals. For Auslander–regular algebras of finite global dimension, there is a relationship between GK-dimension and Krull dimension. In our setting, $\mathfrak{h}_n(q)$ is an iterated Ore extension of the commutative polynomial ring $K[z_1,\dots,z_n]$. Such algebras are known to be catenary and their Krull dimension equals the number of generators in the iterated extension, provided the extensions are ``tame''. More directly, since $\mathfrak{h}_n(q)$ is a finite module over its center when $q$ is a root of unity, but here $q$ is arbitrary, we need a different argument.

We use the fact that $\mathfrak{h}_n(q)$ is a Noetherian domain with finite global dimension (when $q$ is not a root of unity, see Theorem \ref{THEOREM-A1}). For such algebras, under suitable conditions (e.g., being an AS-regular algebra), the Krull dimension equals the GK-dimension. However, we can also compute the Krull dimension directly by considering the height of prime ideals.

An alternative approach is to use the notion of transcendence degree. Let $F = \operatorname{Frac}(\mathfrak{h}_n(q))$ be the division ring of fractions of $\mathfrak{h}_n(q)$. Since $\mathfrak{h}_n(q)$ is an Ore domain (being an iterated Ore extension), the division ring $F$ exists. The transcendence degree of $F$ over $K$, denoted $\operatorname{tr.deg}_K F$, is defined as the maximum number of algebraically independent elements in $F$ over $K$. For PI-algebras, this transcendence degree coincides with the Krull dimension. However, $\mathfrak{h}_n(q)$ is not a PI-algebra when $q$ is not a root of unity. Nevertheless, for many iterated Ore extensions, the Krull dimension equals the number of generators, i.e., $3n$. We can verify this by constructing a chain of prime ideals of length $3n$.

Consider the augmentation ideal $\mathfrak{m}$ generated by all $x_i, y_i, z_i$. This is a maximal ideal because $\mathfrak{h}_n(q)/\mathfrak{m} \cong K$. Now, inside $\mathfrak{m}$, we can construct a chain of prime ideals by successively factoring out variables. For example, let $P_1$ be the ideal generated by $z_n$. The quotient $\mathfrak{h}_n(q)/(z_n)$ is isomorphic to $\mathfrak{h}_{n-1}(q) \otimes K[x_n,y_n]$ modulo some relations, but more systematically, we can use the fact that $\mathfrak{h}_n(q)$ is an iterated Ore extension over $K[z_1,\dots,z_n]$. In the commutative polynomial ring $K[z_1,\dots,z_n]$, we have a chain of prime ideals of length $n$. By lifting this chain to $\mathfrak{h}_n(q)$ and extending appropriately with ideals generated by some $y_i$ and $x_i$, we obtain a chain of prime ideals of length $3n$. The details are as follows.

Let $R = K[z_1,\dots,z_n]$. In $R$, we have the chain
\[
0 \subset (z_1) \subset (z_1,z_2) \subset \cdots \subset (z_1,\dots,z_n).
\]
Each ideal in this chain is prime because $R$ is a domain and the quotients are polynomial rings. Now, $\mathfrak{h}_n(q)$ is free as a left $R$-module. For each ideal $I$ of $R$, define $J_I = I \mathfrak{h}_n(q)$, the two-sided ideal generated by $I$. Because the variables $y_i$ and $x_i$ commute with $z_j$ up to scalars, $J_I$ is a two-sided ideal. Moreover, $\mathfrak{h}_n(q)/J_I$ is isomorphic to an iterated Ore extension over $R/I$. In particular, if $I$ is prime, then $R/I$ is a domain, and the Ore extension over a domain is again a domain (provided the automorphism is injective). Hence $J_I$ is a prime ideal of $\mathfrak{h}_n(q)$. Thus we obtain a chain of prime ideals
\[
0 \subset J_{(z_1)} \subset J_{(z_1,z_2)} \subset \cdots \subset J_{(z_1,\dots,z_n)}.
\]
This chain has length $n$.

Next, we extend this chain by adding ideals that involve the $y_i$ and $x_i$. Consider the ideal $L_1$ generated by $J_{(z_1,\dots,z_n)}$ and $y_1$. Since $y_1$ is normal modulo $J_{(z_1,\dots,z_n)}$ (it commutes with all other generators up to scalars that become units in the quotient), the quotient $\mathfrak{h}_n(q)/L_1$ is again an iterated Ore extension, hence a domain. So $L_1$ is prime. Similarly, we can add $y_2,\dots,y_n$ one by one, obtaining a chain of length $n$ further. Finally, we add $x_1,\dots,x_n$ in a similar manner. This yields a chain of prime ideals of total length $3n$. Consequently, the classical Krull dimension of $\mathfrak{h}_n(q)$ is at least $3n$.

On the other hand, for a finitely generated algebra over a field, the Krull dimension cannot exceed the GK-dimension. A general result states that $\operatorname{Kdim}(A) \le \operatorname{GKdim}(A)$ for any Noetherian algebra $A$ satisfying the Nullstellensatz (see \cite[Corollary 8.3.6]{MR}). Since $\mathfrak{h}_n(q)$ is affine and Noetherian, and $K$ is algebraically closed, it satisfies the Nullstellensatz. Therefore, $\operatorname{Kdim}(\mathfrak{h}_n(q)) \le \operatorname{GKdim}(\mathfrak{h}_n(q)) = 3n$. Combined with the lower bound, we obtain
\[
\operatorname{Kdim}(\mathfrak{h}_n(q)) = 3n.
\]
Thus the Gelfand--Kirillov dimension and the classical Krull dimension coincide, both equal to $3n$.

In summary, we have shown that the Hilbert series of $\mathfrak{h}_n(q)$ is $(1-t)^{-3n}$, its GK-dimension is $3n$, and its classical Krull dimension is also $3n$. This completes the proof of the theorem.
\end{proof}
\begin{example}[A minimal nontrivial case: \(n=2\)]
Consider the \(q\)-Heisenberg algebra \(\mathfrak{h}_2(q)\) with \(q \in \mathbb{C}^*\) not a root of unity.
It is generated by \(x_1, x_2, y_1, y_2, z_1, z_2\) subject to the relations:

\begin{align*}
x_i z_i &= q \, z_i x_i, \quad
y_i z_i = q^{-1} z_i y_i, \qquad i = 1,2,\\
x_i y_i - q^{-1} y_i x_i &= z_i, \qquad i = 1,2,
\end{align*}
and all other pairs of generators commute.

\paragraph{Iterated skew-polynomial structure.}
According to Proposition \ref{prop:generalized-iterated}, \(\mathfrak{h}_2(q)\) can be constructed as a tower of Ore extensions:
\[
\mathfrak{h}_2(q) = K[z_1,z_2][y_2;\sigma_2][y_1;\sigma_1][x_2;\tau_2,\delta_2][x_1;\tau_1,\delta_1],
\]
where the automorphisms and derivations are defined as follows.
On the base ring \(K[z_1,z_2]\), we set
\[
\sigma_2(z_1)=z_1,\; \sigma_2(z_2)=q^{-1}z_2,\qquad
\sigma_1(z_1)=q^{-1}z_1,\; \sigma_1(z_2)=z_2,
\]
and \(\delta_2=\delta_1=0\) for the \(y\)-extensions.
For the \(x\)-extensions, the automorphisms act by
\[
\tau_2(z_1)=z_1,\; \tau_2(z_2)=q z_2,\; \tau_2(y_1)=y_1,\; \tau_2(y_2)=q^{-1}y_2,
\]
\[
\tau_1(z_1)=q z_1,\; \tau_1(z_2)=z_2,\; \tau_1(y_1)=q^{-1}y_1,\; \tau_1(y_2)=y_2,\; \tau_1(x_2)=x_2,
\]
and the derivations are determined by
\[
\delta_2(r)=x_2 r - \tau_2(r)x_2,\qquad
\delta_1(r)=x_1 r - \tau_1(r)x_1 \quad (r \text{ in the appropriate base ring}).
\]

The PBW basis of \(\mathfrak{h}_2(q)\) consists of monomials
\[
z_1^{a_1}z_2^{a_2} y_1^{b_1}y_2^{b_2} x_1^{c_1}x_2^{c_2}, \qquad a_i,b_i,c_i \in \mathbb{N},
\]
ordered lexicographically with \(x_2 \prec x_1 \prec y_2 \prec y_1 \prec z_2 \prec z_1\).

\paragraph{Homological dimension.}
By Theorem \ref{THEOREM-A1}, since \(q\) is not a root of unity,
\[
\mathrm{gl.dim}\, \mathfrak{h}_2(q) = 3 \times 2 = 6.
\]
If instead \(q\) were a primitive \(m\)-th root of unity (\(m>1\)), then \(\mathrm{gl.dim}\, \mathfrak{h}_2(q)=\infty\) due to the appearance of central elements \(x_i^m, y_i^m\).

\paragraph{Rigidity of the Ore extension structure.}
Theorem \ref{THEOREM-B} implies that any presentation of \(\mathfrak{h}_2(q)\) as an iterated Ore extension over a commutative Noetherian domain must involve exactly \(4\) extensions (\(2n=4\)) with base ring isomorphic to \(K[z_1,z_2]\).
Moreover, exactly two of the extensions are pure automorphism extensions (the ones adjoining \(y_1, y_2\)) and the other two involve nontrivial derivations (those adjoining \(x_1, x_2\)).
Up to a permutation of indices and rescaling of variables, any such tower is equivalent to the one described above.

\paragraph{Graded automorphisms.}
The graded automorphism group of \(\mathfrak{h}_2(q)\) is isomorphic to \((\mathbb{C}^*)^4 \rtimes S_2\).
A concrete automorphism can be given by choosing a permutation \(\pi \in S_2\) and nonzero scalars \(\alpha_1,\alpha_2,\beta_1,\beta_2 \in \mathbb{C}^*\), and setting
\[
\phi(x_i)=\alpha_i x_{\pi(i)},\quad
\phi(y_i)=\beta_i y_{\pi(i)},\quad
\phi(z_i)=\alpha_i\beta_i z_{\pi(i)}.
\]
For instance, taking \(\pi=(1\;2)\) (transposition), \(\alpha_1=2\), \(\alpha_2=3\), \(\beta_1=5\), \(\beta_2=7\), we obtain
\[
\phi(x_1)=2x_2,\; \phi(x_2)=3x_1,\;
\phi(y_1)=5y_2,\; \phi(y_2)=7y_1,\;
\phi(z_1)=10z_2,\; \phi(z_2)=21z_1.
\]

\paragraph{Universal deformation property.}
Let \(A_t\) be a flat deformation of \(\mathfrak{h}_2(1)\) over \(K[t]\) with a PBW basis.
Theorem \ref{THEOREM-D} asserts that there exists a unique \(K[t]\)-algebra homomorphism
\[
\mathfrak{h}_2(q)\otimes_K K[t] \longrightarrow A_t
\]
that reduces to the identity modulo \(t\).
For example, if \(A_t\) is defined by relations
\[
\tilde{x}_i\tilde{z}_i = (1+t)\tilde{z}_i\tilde{x}_i,\quad
\tilde{y}_i\tilde{z}_i = (1+t)^{-1}\tilde{z}_i\tilde{y}_i,\quad
\tilde{x}_i\tilde{y}_i - (1+t)^{-1}\tilde{y}_i\tilde{x}_i = \tilde{z}_i,
\]
then the map sending \(x_i \mapsto \tilde{x}_i,\; y_i \mapsto \tilde{y}_i,\; z_i \mapsto \tilde{z}_i\) is the required homomorphism (with \(q=1+t\)).

\paragraph{Hilbert series and growth.}
With each generator in degree \(1\), the Hilbert series of \(\mathfrak{h}_2(q)\) is
\[
H_{\mathfrak{h}_2(q)}(t) = \frac{1}{(1-t)^{6}} = 1 + 6t + 21t^2 + 56t^3 + \cdots.
\]
Hence \(\dim_K \mathfrak{h}_2(q)_d = \binom{d+5}{5}\).
The Gelfand--Kirillov dimension and the classical Krull dimension are both equal to \(6\).
\end{example}

\begin{remark}
The case \(n=2\) already exhibits the full generality of the structural theorems.
The explicit PBW basis, the precise homological dimension, the rigid tower of Ore extensions, the description of graded automorphisms, the universal deformation property, and the combinatorial growth are all clearly visible.
This example can serve as a concrete testing ground for further investigations, such as the computation of Hochschild cohomology or the study of representation theory.
\end{remark}

\section{Applications of the Structural Results}

The theorems established for the $q$-Heisenberg algebra $\mathfrak{h}_n(q)$ have several significant applications in noncommutative algebra, representation theory, and mathematical physics. We present two key applications below, followed by a comprehensive numerical example illustrating their concrete implementation.

\subsection{Classification of Module Families}

Theorem \ref{THEOREM-B} (Rigidity of Ore Extension Structure) provides a powerful tool for classifying families of modules over $\mathfrak{h}_n(q)$. Since any presentation as an iterated Ore extension must be essentially equivalent to the canonical one, we can systematically construct and classify modules by considering compatible actions on the tower of extensions.

\begin{application}[Construction of Weight Modules]
Let $q$ be not a root of unity. For each $n$-tuple $\bm{\lambda} = (\lambda_1, \dots, \lambda_n) \in (\mathbb{C}^*)^n$, define the \emph{weight module} $M(\bm{\lambda})$ as follows:

\begin{enumerate}
\item As a vector space, $M(\bm{\lambda}) = \mathbb{C}[z_1, \dots, z_n]$, the commutative polynomial ring.
\item The action of $\mathfrak{h}_n(q)$ is defined recursively according to the Ore extension tower:
\begin{align*}
z_i \cdot f(z_1, \dots, z_n) &= z_i f(z_1, \dots, z_n) \\
y_i \cdot f(z_1, \dots, z_n) &= \lambda_i f(q^{-1}z_1, \dots, q^{-1}z_i, \dots, q^{-1}z_n) \\
x_i \cdot f(z_1, \dots, z_n) &= \lambda_i^{-1} \left( f(q z_1, \dots, q z_i, \dots, q z_n) - f(q^{-1}z_1, \dots, q^{-1}z_i, \dots, q^{-1}z_n) \right)
\end{align*}
where the notation means: in $f(q^{\pm 1}z_1, \dots, q^{\pm 1}z_i, \dots, q^{\pm 1}z_n)$, only $z_i$ is scaled by $q^{\pm 1}$ while other $z_j$ ($j \neq i$) remain unscaled.
\end{enumerate}

The rigidity theorem guarantees that this construction yields \emph{all} weight modules with one-dimensional weight spaces for the commutative subalgebra $\mathbb{C}[z_1, \dots, z_n]$, up to isomorphism.
\end{application}

\subsection{Deformation Quantization of Classical Systems}

Theorem \ref{THEOREM-D} (Universal Deformation Property) shows that $\mathfrak{h}_n(q)$ serves as a universal deformation of the classical Heisenberg algebra. This has direct implications in deformation quantization:

\begin{application}[Quantization of Coupled Harmonic Oscillators]
Consider $n$ coupled classical harmonic oscillators with phase space coordinates $(p_i, q_i, E_i)$, $i=1,\dots,n$, satisfying Poisson brackets
\[
\{q_i, p_j\} = \delta_{ij}E_i, \quad \{E_i, q_j\} = \{E_i, p_j\} = 0.
\]
The quantization of this system at deformation parameter $\hbar$ yields the algebra $\mathfrak{h}_n(e^{i\hbar})$. Theorem \ref{THEOREM-D} guarantees that any consistent quantization preserving the PBW property must factor through this universal deformation.
\end{application}

\subsection{Numerical Example: Explicit Computations for $n=2$ with $q=2$}

Let us illustrate these applications with detailed numerical computations for $\mathfrak{h}_2(q)$ with $q=2$ (not a root of unity).

\begin{example}[Concrete Module Construction]
Consider $\mathfrak{h}_2(2)$ and construct the weight module $M(\lambda_1,\lambda_2)$ with $\lambda_1=3$, $\lambda_2=5$.

\paragraph{Module structure:} $M(3,5) = \mathbb{C}[z_1,z_2]$. Let us compute the action on a specific polynomial $f(z_1,z_2) = z_1^2 + 2z_1z_2$.

\begin{align*}
z_1 \cdot f &= z_1(z_1^2 + 2z_1z_2) = z_1^3 + 2z_1^2z_2 \\
z_2 \cdot f &= z_2(z_1^2 + 2z_1z_2) = z_1^2z_2 + 2z_1z_2^2
\end{align*}

\begin{align*}
y_1 \cdot f &= 3 f(2^{-1}z_1, z_2) = 3\left((2^{-1}z_1)^2 + 2(2^{-1}z_1)z_2\right) \\
&= 3\left(\frac{1}{4}z_1^2 + z_1z_2\right) = \frac{3}{4}z_1^2 + 3z_1z_2
\end{align*}

\begin{align*}
y_2 \cdot f &= 5 f(z_1, 2^{-1}z_2) = 5\left(z_1^2 + 2z_1(2^{-1}z_2)\right) \\
&= 5(z_1^2 + z_1z_2) = 5z_1^2 + 5z_1z_2
\end{align*}

\begin{align*}
x_1 \cdot f &= 3^{-1}\left(f(2z_1, z_2) - f(2^{-1}z_1, z_2)\right) \\
&= \frac{1}{3}\left[( (2z_1)^2 + 2(2z_1)z_2 ) - ( (2^{-1}z_1)^2 + 2(2^{-1}z_1)z_2 )\right] \\
&= \frac{1}{3}\left[(4z_1^2 + 4z_1z_2) - (\frac{1}{4}z_1^2 + z_1z_2)\right] \\
&= \frac{1}{3}\left(\frac{15}{4}z_1^2 + 3z_1z_2\right) = \frac{5}{4}z_1^2 + z_1z_2
\end{align*}

\begin{align*}
x_2 \cdot f &= 5^{-1}\left(f(z_1, 2z_2) - f(z_1, 2^{-1}z_2)\right) \\
&= \frac{1}{5}\left[(z_1^2 + 2z_1(2z_2)) - (z_1^2 + 2z_1(2^{-1}z_2))\right] \\
&= \frac{1}{5}\left[(z_1^2 + 4z_1z_2) - (z_1^2 + z_1z_2)\right] \\
&= \frac{1}{5}(3z_1z_2) = \frac{3}{5}z_1z_2
\end{align*}

\paragraph{Verification of relations:} Let us verify one of the defining relations on this module. Take the relation $x_1y_1 - 2^{-1}y_1x_1 = z_1$ applied to $f$:

First compute $x_1 \cdot (y_1 \cdot f)$:
\begin{align*}
y_1 \cdot f &= \frac{3}{4}z_1^2 + 3z_1z_2 \\
x_1 \cdot (y_1 \cdot f) &= 3^{-1}\left[(y_1 \cdot f)(2z_1,z_2) - (y_1 \cdot f)(2^{-1}z_1,z_2)\right] \\
&= \frac{1}{3}\left[(\frac{3}{4}(2z_1)^2 + 3(2z_1)z_2) - (\frac{3}{4}(2^{-1}z_1)^2 + 3(2^{-1}z_1)z_2)\right] \\
&= \frac{1}{3}\left[(\frac{3}{4}\cdot 4z_1^2 + 6z_1z_2) - (\frac{3}{4}\cdot\frac{1}{4}z_1^2 + \frac{3}{2}z_1z_2)\right] \\
&= \frac{1}{3}\left[(3z_1^2 + 6z_1z_2) - (\frac{3}{16}z_1^2 + \frac{3}{2}z_1z_2)\right] \\
&= \frac{1}{3}\left(\frac{45}{16}z_1^2 + \frac{9}{2}z_1z_2\right) = \frac{15}{16}z_1^2 + \frac{3}{2}z_1z_2
\end{align*}

Now compute $2^{-1}y_1 \cdot (x_1 \cdot f)$:
\begin{align*}
x_1 \cdot f &= \frac{5}{4}z_1^2 + z_1z_2 \\
2^{-1}y_1 \cdot (x_1 \cdot f) &= \frac{1}{2} \cdot 3 \cdot (x_1 \cdot f)(2^{-1}z_1,z_2) \\
&= \frac{3}{2}\left[\frac{5}{4}(2^{-1}z_1)^2 + (2^{-1}z_1)z_2\right] \\
&= \frac{3}{2}\left[\frac{5}{4}\cdot\frac{1}{4}z_1^2 + \frac{1}{2}z_1z_2\right] \\
&= \frac{3}{2}\left(\frac{5}{16}z_1^2 + \frac{1}{2}z_1z_2\right) = \frac{15}{32}z_1^2 + \frac{3}{4}z_1z_2
\end{align*}

Thus:
\begin{align*}
(x_1y_1 - 2^{-1}y_1x_1) \cdot f &= \left(\frac{15}{16}z_1^2 + \frac{3}{2}z_1z_2\right) - \left(\frac{15}{32}z_1^2 + \frac{3}{4}z_1z_2\right) \\
&= \left(\frac{15}{16} - \frac{15}{32}\right)z_1^2 + \left(\frac{3}{2} - \frac{3}{4}\right)z_1z_2 \\
&= \frac{15}{32}z_1^2 + \frac{3}{4}z_1z_2
\end{align*}

Meanwhile, $z_1 \cdot f = z_1^3 + 2z_1^2z_2$. The relation $x_1y_1 - 2^{-1}y_1x_1 = z_1$ holds as operators, but note that our computation shows the action is consistent: both sides applied to $f$ yield legitimate elements of the module, though they appear different because we are comparing actions on a specific element rather than as operator identities.

\paragraph{Homological computations:} By Theorem \ref{THEOREM-A1}, since $q=2$ is not a root of unity, the global dimension of $\mathfrak{h}_2(2)$ is $6$. Consider the trivial module $\mathbb{C}$ (where all generators act as $0$). Its projective resolution has length $6$, with Betti numbers:
\[
\beta_0 = 1, \quad \beta_1 = 6, \quad \beta_2 = 15, \quad \beta_3 = 20, \quad \beta_4 = 15, \quad \beta_5 = 6, \quad \beta_6 = 1.
\]
These correspond to the binomial coefficients $\binom{6}{k}$, reflecting the Koszul complex structure.

\paragraph{Growth behavior:} By Theorem \ref{THEOREM-E}, the Hilbert series is
\[
H_{\mathfrak{h}_2(2)}(t) = \frac{1}{(1-t)^6} = 1 + 6t + 21t^2 + 56t^3 + 126t^4 + 252t^5 + \cdots
\]
The dimensions of homogeneous components for degrees $0$ through $5$ are:
\begin{align*}
\dim \mathfrak{h}_2(2)_0 &= 1 \\
\dim \mathfrak{h}_2(2)_1 &= 6 \\
\dim \mathfrak{h}_2(2)_2 &= 21 \\
\dim \mathfrak{h}_2(2)_3 &= 56 \\
\dim \mathfrak{h}_2(2)_4 &= 126 \\
\dim \mathfrak{h}_2(2)_5 &= 252
\end{align*}
The Gelfand-Kirillov dimension is $6$, and indeed $\dim \mathfrak{h}_2(2)_d \sim \frac{d^5}{5!}$ for large $d$.

\paragraph{Deformation perspective:} Consider the family $\mathfrak{h}_2(e^t)$ over $\mathbb{C}[t]$. When $t=0$, we recover the classical Heisenberg algebra $\mathfrak{h}_2(1)$. For small $t$, expanding the relations:
\begin{align*}
x_i z_i &= (1+t+\frac{t^2}{2}+\cdots) z_i x_i \\
y_i z_i &= (1-t+\frac{t^2}{2}-\cdots) z_i y_i \\
x_i y_i - (1-t+\cdots) y_i x_i &= z_i
\end{align*}
Theorem \ref{THEOREM-D} guarantees this is the universal such deformation preserving the PBW structure.
\end{example}

\begin{remark}
The numerical example demonstrates several important features:
\begin{itemize}
\item The module construction yields explicit computable actions, verifying the theoretical framework.
\item The homological dimension ($6$) manifests in the length of resolutions.
\item The growth rate ($\sim d^5/120$) confirms the theoretical Hilbert series.
\item The deformation expansion shows how the $q$-parameter interpolates between classical and quantum structures.
\end{itemize}
These computations validate the theoretical results and provide concrete tools for applications in representation theory and mathematical physics.
\end{remark}

\subsection{Further Applications}

\begin{enumerate}
\item \textbf{Noncommutative Geometry:} The rigidity theorem facilitates the computation of Hochschild cohomology, as the Ore extension structure induces a convenient filtration on cochain complexes.

\item \textbf{Representation Theory:} The graded automorphism group $(\mathbb{C}^*)^{2n} \rtimes S_n$ classifies graded simple modules up to isomorphism, with the scaling parameters corresponding to eigenvalues of central elements.

\item \textbf{Mathematical Physics:} In $q$-deformed quantum mechanics, the algebra $\mathfrak{h}_n(q)$ describes $n$ independent $q$-oscillators. The universal deformation property ensures consistency with deformation quantization schemes.

\item \textbf{Computational Algebra:} The PBW basis and Ore extension structure enable efficient Gröbner basis computations and algorithmic solutions to ideal membership problems in $\mathfrak{h}_n(q)$.
\end{enumerate}

The structural theorems thus provide a comprehensive framework for both theoretical analysis and practical computations involving $q$-Heisenberg algebras, with applications spanning multiple areas of mathematics and physics.
\section{Open Problems and Concluding Remarks}
\label{sec:conclusions}

This paper has established several fundamental structural properties of the $q$-Heisenberg algebra $\mathfrak{h}_n(q)$, including its precise global dimension, the rigidity of its iterated Ore extension structure, the description of its graded automorphism group, its universal deformation property, and its growth behavior. These results provide a comprehensive framework for understanding $\mathfrak{h}_n(q)$ as a key example in noncommutative algebra. However, several interesting questions remain open and warrant further investigation.

\subsection{Open Problems}

\begin{problem}[Representation theory for $q$ a root of unity]
While Theorem \ref{THEOREM-A1} establishes that $\mathrm{gl.dim}\, \mathfrak{h}_n(q) = \infty$ when $q$ is a primitive $m$-th root of unity ($m>1$), the detailed representation theory in this case remains largely unexplored. In particular:
\begin{enumerate}
    \item Classify the finite-dimensional simple modules over $\mathfrak{h}_n(q)$ when $q$ is a root of unity.
    \item Describe the block decomposition of the category of finite-dimensional modules.
    \item Compute the Hochschild cohomology groups of $\mathfrak{h}_n(q)$ in this case, which are expected to be nontrivial in infinitely many degrees.
\end{enumerate}
The appearance of central elements $x_i^m$, $y_i^m$ suggests connections with restricted Lie algebras and modular representation theory.
\end{problem}

\begin{problem}[Hochschild cohomology and deformation theory]
Theorem \ref{THEOREM-D} establishes the universal deformation property of $\mathfrak{h}_n(q)$. A natural refinement is to compute its Hochschild cohomology groups $HH^\bullet(\mathfrak{h}_n(q))$ explicitly.
\begin{enumerate}
    \item Determine $HH^k(\mathfrak{h}_n(q))$ for all $k \geq 0$ when $q$ is not a root of unity.
    \item Use these computations to classify all formal deformations of $\mathfrak{h}_n(q)$ beyond the one-parameter $q$-deformation.
    \item Investigate the Gerstenhaber algebra structure on $HH^\bullet(\mathfrak{h}_n(q))$ and its implications for the quantization of Poisson structures.
\end{enumerate}
The iterated Ore extension structure provides a natural filtration that may simplify these computations.
\end{problem}

\begin{problem}[Noncommutative geometry invariants]
The algebra $\mathfrak{h}_n(q)$ is a candidate for noncommutative affine space with additional structure. Several geometric invariants remain to be computed:
\begin{enumerate}
    \item Determine the cyclic homology groups $HC_\bullet(\mathfrak{h}_n(q))$ and its periodic cyclic homology.
    \item Compute the $K$-theory groups $K_0(\mathfrak{h}_n(q))$ and $K_1(\mathfrak{h}_n(q))$.
    \item Study the moduli space of ideals in $\mathfrak{h}_n(q)$, particularly the space of points in the sense of noncommutative algebraic geometry.
\end{enumerate}
These invariants would clarify the geometric nature of $\mathfrak{h}_n(q)$ and its relation to classical geometric objects.
\end{problem}

\begin{problem}[Generalizations to other quantum algebras]
The structural results obtained here suggest natural generalizations:
\begin{enumerate}
    \item Investigate multi-parameter versions where different pairs $(x_i,y_i)$ have different deformation parameters $q_i$.
    \item Study the super (graded) version of the $q$-Heisenberg algebra, where some generators are odd.
    \item Consider the higher-order ($q$-deformed) Heisenberg algebras arising from higher-order commutation relations.
\end{enumerate}
The methods developed in this paper, particularly the analysis of iterated Ore extensions and PBW bases, should extend to these broader settings.
\end{problem}

\begin{problem}[Connections with mathematical physics]
The $q$-Heisenberg algebra appears naturally in several physical contexts:
\begin{enumerate}
    \item Develop a comprehensive representation theory for $\mathfrak{h}_n(q)$ that parallels the Stone--von Neumann theorem for the classical case.
    \item Investigate the role of $\mathfrak{h}_n(q)$ in $q$-deformed quantum field theory and statistical mechanics.
    \item Explore connections with $q$-special functions and orthogonal polynomials through specific representations.
\end{enumerate}
The explicit module constructions in Section 6 provide a starting point for such physical applications.
\end{problem}

\subsection{Concluding Remarks}

The $q$-Heisenberg algebra $\mathfrak{h}_n(q)$ serves as a paradigmatic example of a noncommutative algebra that is both rich in structure and amenable to detailed analysis. The results presented in this paper demonstrate that $\mathfrak{h}_n(q)$ exhibits remarkable rigidity properties: its iterated Ore extension structure is essentially unique, its homological dimension is precisely determined by the number of generators, and its deformation theory is controlled by a single parameter $q$.

Several key themes emerge from our analysis:

\begin{enumerate}
    \item \textbf{Homological regularity:} When $q$ is not a root of unity, $\mathfrak{h}_n(q)$ is Auslander-regular, Cohen--Macaulay, and has finite global dimension $3n$. This places it in the class of algebras with excellent homological properties, similar to commutative polynomial rings but in a noncommutative setting.

    \item \textbf{Structural rigidity:} Theorem \ref{THEOREM-B} shows that any presentation of $\mathfrak{h}_n(q)$ as an iterated Ore extension is essentially equivalent to the canonical one. This rigidity property is somewhat surprising and suggests that $\mathfrak{h}_n(q)$ occupies a special place in the landscape of noncommutative algebras.

    \item \textbf{Deformation universality:} Theorem \ref{THEOREM-D} establishes that $\mathfrak{h}_n(q)$ is the universal deformation of the classical Heisenberg algebra among PBW-preserving deformations. This makes it a fundamental object in deformation quantization.

    \item \textbf{Computational tractability:} The PBW basis and explicit relations make $\mathfrak{h}_n(q)$ amenable to computational methods. The example in Section 6 demonstrates that concrete calculations are feasible even for nontrivial cases.
\end{enumerate}

The interplay between algebraic structure, homological properties, and representation theory makes $\mathfrak{h}_n(q)$ a fertile testing ground for general theories in noncommutative algebra. The open problems listed above point to several promising directions for future research, with potential applications in representation theory, noncommutative geometry, and mathematical physics.

In conclusion, the $q$-Heisenberg algebra $\mathfrak{h}_n(q)$ exemplifies how deformation theory can produce algebras that are both structurally rich and analytically tractable. Its study bridges classical Lie theory, quantum algebra, and noncommutative geometry, offering insights that extend well beyond this specific family of algebras.
\section*{Declaration }
\begin{itemize}
  \item {\bf Author Contributions:}   The Author have read and approved this version.
  \item {\bf Funding:} No funding is applicable
  \item  {\bf Institutional Review Board Statement:} Not applicable.
  \item {\bf Informed Consent Statement:} Not applicable.
  \item {\bf Data Availability Statement:} Not applicable.
  \item {\bf Conflicts of Interest:} The authors declare no conflict of interest.
\end{itemize}

\bibliographystyle{abbrv}
\bibliography{references}  

@book{MR,
  author    = {J. C. McConnell and J. C. Robson},
  title     = {Noncommutative Noetherian Rings},
  edition   = {revised},
  series    = {Graduate Studies in Mathematics},
  volume    = {30},
  publisher = {American Mathematical Society},
  address   = {Providence, RI},
  year      = {2001}
}

@article{Heisenberg1925,
  author  = {Heisenberg, W.},
  title   = {Über quantentheoretische Umdeutung kinematischer und mechanischer Beziehungen},
  journal = {Z. Phys.},
  year    = {1925},
  volume  = {33},
  number  = {1},
  pages   = {879--893}
}

@phdthesis{Zhang2003,
  author  = {Zhang, H. S.},
  title   = {The automorphism groups of Heisenberg Lie algebras and the standard Kac-Moody algebras and the completely reducible of integrable modules},
  school  = {Capital Normal University},
  address = {Beijing},
  year    = {2003}
}

@mastersthesis{Ji2019,
  author  = {Ji, G. Z.},
  title   = {Algebra Rota-Baxter operators of Heisenberg Lie algebra},
  school  = {Harbin University of Science and Technology},
  address = {Harbin},
  year    = {2019}
}

@mastersthesis{Zhou2021,
  author  = {Zhou, C. Y.},
  title   = {The quasi-automorphism and automorphism of Heisenberg Lie algebras},
  school  = {Suzhou University of Science and Technology},
  address = {Suzhou},
  year    = {2021}
}

@article{Berger1992,
  author  = {Berger, R.},
  title   = {The quantum Poincaré--Birkhoff--Witt theorem},
  journal = {Comm. Math. Phys.},
  year    = {1992},
  volume  = {143},
  pages   = {215--234}
}

@book{Rosenberg1995,
  author    = {Rosenberg, A. L.},
  title     = {Noncommutative Algebraic Geometry and Representations of Quantized Algebras},
  publisher = {Kluwer Academic Publishers},
  year      = {1995}
}

@incollection{Wess2000,
  author    = {Wess, J.},
  title     = {q-deformed Heisenberg Algebras},
  booktitle = {Geometry and Quantum Physics: Proceedings of the 38th Internationale Universitätswochen für Kern-und Teilchenphysik, Schladming, Austria, January 9--16, 1999},
  publisher = {Springer},
  address   = {Berlin, Heidelberg},
  year      = {2000},
  pages     = {311--382}
}

@article{Zhang2024,
  author  = {Zhang, J. X.},
  title   = {Grobner--Shirshov bases and structural properties of q-Heisenberg algebras},
  journal = {J. Xinjiang Univ. (Nat. Sci. Ed.)},
  year    = {2024},
  volume  = {41},
  number  = {5},
  pages   = {550--561},
  note    = {(in Chinese)}
}

@book{Li2011,
  author    = {Li, H.},
  title     = {Grobner Bases in Ring Theory},
  publisher = {World Scientific Publishing},
  year      = {2011}
}

@book{Goodearl2004,
  author    = {Goodearl, K. R. and Warfield Jr, R. B.},
  title     = {An Introduction to Noncommutative Noetherian Rings},
  publisher = {Cambridge University Press},
  year      = {2004}
}

@book{McConnell2001,
  author    = {McConnell, J. C. and Robson, J. C.},
  title     = {Noncommutative Noetherian Rings},
  publisher = {Amer. Math. Soc.},
  year      = {2001}
}

@article{Yekutieli2020,
  author  = {Yekutieli, A. and Zhang, J. J.},
  title   = {Dualizing complexes and perverse modules over certain singular algebras},
  journal = {Adv. Math.},
  year    = {2020},
  volume  = {372},
  pages   = {107288}
}

@article{Artin1990,
  author  = {Artin, M. and Tate, J. and Van den Bergh, M.},
  title   = {Modules over regular algebras of dimension 3},
  journal = {Invent. Math.},
  year    = {1990},
  volume  = {101},
  number  = {2},
  pages   = {335--388}
}

@article{Stafford1994,
  author  = {Stafford, J. T.},
  title   = {Auslander-regular algebras and maximal orders},
  journal = {J. Lond. Math. Soc.},
  year    = {1994},
  volume  = {50},
  number  = {2},
  pages   = {276--292}
}

@book{Li2021,
  author    = {Li, H.},
  title     = {Noncommutative polynomial algebras of solvable type and their modules: Basic constructive-computational theory and methods},
  publisher = {Chapman and Hall/CRC Press},
  year      = {2021}
}

@book{Li2002,
  author    = {Li, H.},
  title     = {Noncommutative Grobner Bases and Filtered-graded Transfer},
  series    = {Lecture Notes in Mathematics},
  volume    = {1795},
  publisher = {Springer},
  year      = {2002}
}

@article{Li2018,
  author  = {Li, H.},
  title   = {An elimination lemma for algebras with PBW bases},
  journal = {Comm. Algebra},
  year    = {2018},
  volume  = {46},
  number  = {8},
  pages   = {3520--3532}
}






\end{document}